\documentclass[a4paper,10pt,final]{amsart}
\usepackage{amsmath}
\usepackage{amsfonts}
\usepackage{amssymb}
\usepackage{amsthm}
\usepackage{graphicx}
\usepackage{amsmath}
\usepackage{psfrag}
\usepackage{cite}
\usepackage[dvips]{hyperref}
\usepackage{showkeys}

\newtheorem{theorem}{Theorem}[section]
\newtheorem{lemma}[theorem]{Lemma}
\newtheorem{remark}[theorem]{Remark}

\newtheorem{corollary}[theorem]{Corollary}
\newtheorem{definition}[theorem]{Definition}

  \newcommand{\DivSymb}{\mathsf{div}}

  \newcommand{\Div}{\DivSymb}

  \newcommand{\fbf}{{\mathbf{f}}}

  \newcommand{\nbf}{{\mathbf{n}}}

  \newcommand{\ubf}{{\mathbf{u}}}
  \newcommand{\vbf}{{\mathbf{v}}}

  \newcommand{\Cbf}{{\mathbf{C}}}

  \newcommand{\Mbf}{{\mathbf{M}}}

  \newcommand{\Rbf}{{\mathbf{R}}}

  \renewcommand{\theta}{\vartheta}
  \renewcommand{\phi}{\varphi}
  \renewcommand{\epsilon}{\varepsilon}

\newcommand{\curl}{\operatorname{\bf curl}}
\newcommand{\grad}{\operatorname{\bf grad}}

\title[Discrete convection-diffusion for differential forms]{Eulerian and
  Semi-Lagrangian Methods for Convection-Diffusion for Differential Forms}

\author{Holger Heumann and Ralf Hiptmair}
\thanks{Authors' affiliation: SAM, ETH Z\"urich, CH-8092 Z\"urich, \{hiptmair,heumann\}@sam.math.ethz.ch}

\date{\today, submitted to special issue of Disc. Cont. Dyn. Sys.}

\begin{document}
\maketitle
\begin{abstract}
  We consider generalized linear transient convection-diffusion problems for
  differential forms on bounded domains in $\mathbb{R}^{n}$. These involve Lie
  derivatives with respect to a prescribed smooth vector field. We construct both new
  Eulerian and semi-Lagrangian approaches to the discretization of the Lie
  derivatives in the context of a Galerkin approximation based on discrete
  differential forms. Details of implementation are discussed as well as an
  application to the discretization of eddy current equations in moving media.
\end{abstract}

\section{Introduction}
\label{sec:introduction}

Recall the classical linear transient 2nd-order convection-diffusion problem for an unknown
scalar function $u=u(x,t)$ on a bounded domain
$\Omega\subset\mathbb{R}^{n}$:
\begin{gather}
  \label{eq:cd}
  \begin{array}[c]{rcll}
    \partial_{t}u - \varepsilon \Delta u + {\boldsymbol \beta}\cdot \operatorname{\bf grad} u &=& 
    f & \text{in }\Omega\;,\\
    u &=& 0 & \text{on }\partial\Omega\;,\\
    u(0) &=& u_{0}\;.
  \end{array}
\end{gather}
Here, and in the remainder of the paper,
$\beta:\overline{\Omega}\mapsto\mathbb{R}^{n}$ stands for a smooth vector field. For
$0<\varepsilon\ll 1$ we encounter a singularly perturbed boundary value problem,
whose \emph{stable} discretization has attracted immense attention in numerical
analysis, see \cite{RST96} and the many references cited therein. 

The boundary value problem \eqref{eq:cd} turns out to be a member of a larger family
of 2nd-order boundary value problems (BVP), which can conveniently be described using the
\emph{calculus of differential forms}. It permits us to state the generalized
transient 2nd-order convection-diffusion problems as
\begin{gather}
  \label{eq:cddf}
  \begin{array}[c]{rcll}
    \ast \partial_t \omega(t) - \varepsilon (-1)^{l} 
    d\ast d \omega(t) +  \ast
    L_{\boldsymbol \beta} \omega (t) &=& \varphi \quad &
    \text{in }\Omega\subset\mathbb{R}^{n}\;,\\
    \imath^{\ast}\omega &=& 0 & \text{on }\partial\Omega\;,\\
    \omega(0) &=& \omega_{0}\;.
  \end{array}
\end{gather}

These are BVPs for an unknown time dependent $l$-form $\omega(t)$, $0\leq l\leq n$,
on the domain $\Omega$. The symbol $\ast$ stands for a Hodge operator mapping an
$l$-form to an $n-l$-form, and $d$ denotes the exterior derivative. Dirichlet
boundary conditions are imposed via the trace $\imath^{\ast}\omega$ of the $l$-form
$\omega$. We refer to \cite[Sect.~2]{HIP02} and \cite[Sect.~2]{AFW06} for more
details and an introduction to the calculus of differential forms in the spirit of
this article. By $L_{\boldsymbol \beta}$ we denote the Lie derivative of $\omega$ for
the given velocity field $\boldsymbol \beta$ \cite[Ch.~5]{AMR83}. It provides the
convective part of the differential operator in \eqref{eq:cddf}. Details will be
given in Sect.~\ref{sec:extr-contr-lie} below.

In two and three dimensions differential forms can be modelled by means of functions and
vector fields through so-called vector proxies, see \cite[Sect.~7]{BOS05c},
\cite[Table~2.1]{AFW06}, and \cite[Table~2.1]{HIP02}\footnote{Occasionally we will use
the operator $\operatorname{v.p.}$ to indicate that a form is mapped to its
corresponding Euclidean vector proxy}. Thus, for $n=3$, in the case of
$\ast$ induced by the Euclidean metric on $\mathbb{R}^{3}$, the following vector
analytic incarnations of \eqref{eq:cddf} are obtained. The solution $l$-form $\omega$ is
replaced with an unknown vector field $\ubf$ or function $u$, \textit{cf.}
\cite[Sect.~2]{HIP99c}.
\begin{itemize}
\item For $l=0$ we recover \eqref{eq:cd}. 
\item In the case $l=1$ we arrive at a convection-diffusion problem 
  \begin{gather}
    \label{eq:cd1}
    \begin{array}[c]{rcll}
      \partial_{t}\ubf + \varepsilon\curl\curl\ubf - {\boldsymbol
      \beta}\times \curl \ubf + 
      \operatorname{\bf grad}(\ubf\cdot{\boldsymbol \beta}) &=& \fbf & \text{in
      }\Omega\;,\\
      \ubf\times\nbf &=& 0 & \text{on }\partial\Omega\;,\\
      \ubf(0) &=& \ubf_{0}\;.
    \end{array}
  \end{gather}
\item The corresponding boundary value problems for vector proxies of 2-forms read
   \begin{gather}
    \label{eq:cd2}
    \begin{array}[c]{rcll}
      \partial_{t}\ubf -\varepsilon\grad\Div\ubf + {\boldsymbol \beta}(\Div\ubf) + 
      \operatorname{\bf curl}(\ubf\times{\boldsymbol \beta}) &=& \fbf & \text{in
      }\Omega\;,\\
      \ubf\cdot\nbf &=& 0 & \text{on }\partial\Omega\;,\\
      \ubf(0) &=& \ubf_{0}\;.
    \end{array}
  \end{gather}
\item For $l=3$ the diffusion term vanishes and we end up with pure convection
  problems for a scalar density $u$
   \begin{gather}
    \label{eq:cd3}
    \begin{array}[c]{rcll}
      \partial_{t}u + \Div({\boldsymbol \beta}u)  &=& f & \text{in
      }\Omega\;,\\
      \ubf(\cdot,0) &=& \ubf_{0}\;,
    \end{array}
  \end{gather}
  where we assume ${\boldsymbol \beta}\cdot\nbf=0$ on $\partial\Omega$ 
\end{itemize}

Though equivalent to \eqref{eq:cddf} the vector proxy formulations very much conceal
the common structure inherent in \eqref{eq:cd1}-\eqref{eq:cd2}. Thus, in this
article, we consistently adopt the perspective of differential forms. Thinking in
terms of co-ordinate free differential forms offers considerable benefits as regards
the construction of structure preserving spatial discretizations. By now this
is widely appreciated for boundary value problems for $d\ast d \omega$, where 
the so-called discrete exterior calculus \cite{HIR03,DHM03,AFW09}, or, equivalently,
the mimetic finite difference approach \cite{HYS01,HYS99a,HYS97,BBL08}, or 
discrete Hodge-operators \cite{HIP99c,BOK99a} have shed new light on existing
discretizations and paved the way for new numerical methods. 

All these methods have in common that $l$-forms are approximated by $l$-co-chains
on generalized triangulations of the computational domain $\Omega$. This preserves
\mbox{deRham} co-homology and, thus, plenty of algebraic properties enjoyed by the 
solutions of the continuous BVP carry over to the discrete setting, see, e.g.,
\cite[Sect.~2.2]{HIP02}. Moreover, simplicial and tensor product triangulations
allow the extension of co-chains to discrete differential forms, which furnishes
structure preserving finite element methods. Also in this paper the focus
will be on Galerkin discretization by means of discrete differential forms
and how it can be used to deal with Lie derivatives.

In light of the success of discrete differential forms, it seems worthwhile
exploring their use for the more general equations \eqref{eq:cddf}. This was
pioneered in \cite{BOS04a,MMP08} and this article continues and extends these
considerations. Revealing structure is not our only motivation: as already pointed
out the discretization of \eqref{eq:cd} has been the object of intense research, but
\eqref{eq:cd1} and \eqref{eq:cd2} have received scant attention, though,
\eqref{eq:cd1} is relevant for numerical modelling, e.g., in magnetohydrodynamics.
Apart from standard Galerkin finite element methods and \cite{BOS04a,MMP08}, we would
like to mention \cite{RWW07,DES01} as attempts to devise a meaningful discretization
of the transport term in \eqref{eq:cd1}. The first paper relies on a Lagrangian
approach along with a divergence preserving projection, which can be interpreted as
an interpolation onto co-chains. The second article employs edge degrees of freedom
(1-co-chains) combined with a special transport scheme. Both refrain from utilizing
exterior calculus.

We stress that this article is confined to the derivation and formulation of fully
discrete schemes for \eqref{eq:cddf} on a single \emph{fixed} triangulation of
$\Omega$. Principal attention is paid to the treatment of the convective terms.  The
discussion centers on structural and algorithmic aspects, whereas rigorous numerical
analysis of stability and convergence analysis will be addressed in some forthcoming
work, because it relies on certain technical results that must be established for
each degree $l$ of the form $\omega$ separately. 

The plan of the paper is as follows: We first review the definition of Lie
derivatives and material derivatives for differential forms and give a short summary
on discrete forms. Next, we present a semi-Lagrangian approximation procedure for
material derivatives of arbitrary discrete $l$-forms.  Afterwards we present an
Eulerian treatment of convection of discrete forms. These two sections focus on the
formulation of fully discrete schemes. They are subsequently elaborated for the eddy
current model in moving media as an example for non-scalar convection-diffusion, see
also \cite{HHX09}.  Finally we present numerical experiments for convection of
$1$-forms in $2$ dimensions. They give empiric hints on convergence and stability.

\section{Preliminaries}
\label{sec:preliminaries}

\subsection{Lie derivatives and material derivatives of forms}
\label{sec:extr-contr-lie}

We introduce a space-time domain $Q_T:=\Omega \times [t_0,t_1]$ where $\Omega \subset
\mathbb R^n$ is a bounded curvilinear polyhedron and $[t_0,t_1] \subset
\mathbb R$. Given a vector field $\boldsymbol
\beta: Q_T \mapsto \mathbb R^n$ and initial values $(x,s)\in \bar{Q}_T$ the
associated flow $X(x,s,t): \bar{\Omega} \times [t_0,t_1]^2 \mapsto \mathbb R^n$ is
defined through
\begin{equation}
    \begin{aligned}
	\nonumber \frac{d}{dt}X(x,s,t)&=\boldsymbol \beta(X(x,s,t),t), \quad t
	\in
	[t_0,t_1] \\
	X(x,s,s)&=x\;.
	\label{eqn:charac_ode}
    \end{aligned}
\end{equation}
For fixed $(x,s)$ the solution $X(x,s,\cdot)$ is called the characteristic curve
through $(x,s)$. A unique solution of problem (\ref{eqn:charac_ode}) exists, whenever
$\boldsymbol \beta$ is continuous in $\bar{Q}_T$ and Lipschitz continuous in
$\bar{\Omega}$ for fixed $t \in [t_0,t_1] $. In this case
\begin{eqnarray}
    X (X(x,s,t),t,s)=x, \quad \forall x \in \Omega 
    \label{eqn:charac_inv1}
\end{eqnarray}
hence $X(\cdot,t,s)$ is the inverse of $X(\cdot,s,t)$: 
\begin{eqnarray}
    X(\cdot,t,s)\circ X(\cdot,s,t) = id\,.   
    \label{eqn:charac_inv2}
\end{eqnarray}
For convenience we abbreviate $X_{s,t}(\cdot)=X(\cdot,s,t)$. It follows
directly from \eqref{eqn:charac_ode} that the Jacobian $DX_{s,t}(x)$ solves:
\begin{align}
    \label{eq:Vcharac_ODE}
    \frac{d}{dt} DX_{s,t}(x)=D\boldsymbol \beta (X_{s,t}(x),t)DX_{s,t}(x)\,.
\end{align}
We write $\mathcal{F}^{l}(\Omega)$ for the space of differential $l$-forms on
$\Omega$ (cf.  \cite[Chapter V, Section 3]{Lang}). Differential $l$-forms can be seen
as continuous additive mappings from the set $\mathcal S^l(\Omega)$ of compact
oriented, piecewise smooth, $l$-dimensional sub-manifolds of $\Omega$ into the real
numbers\footnote{This notion of differential forms can be rigorously defined with
  tools from geometric measure theory \cite{ROS03,HAR05}. Here we forgo these
  technical aspects and appeal to intuition, \textit{cf.} \cite[Sect.~2.1]{HIP02}}.
The spaces $\mathcal F^{l}(\Omega)$ can be equipped with the $L^{2}$-norm $\|\omega
\|_0^2:= \int_{\Omega} \ast \omega \wedge \omega$, see \cite[Sect.~2.2]{AFW06}. For
Euclidean Hodge operator this agrees with the usual $L^2$-norm of the vector proxies.
Although it is very common to regard the mapping $\omega : \mathcal S^l(\Omega)
\mapsto \mathbb C, \omega \in \mathcal F^l(\Omega)$ as integration
we adopt the notation $<\omega,M_l>$ of a duality pairing instead of
$\int\nolimits_{M_l} \omega, M_l \in \mathcal S^l(\Omega)$. Only for $l=n$ we will
use the integral symbol.  

\begin{figure}
    \begin{center}
	\psfrag{beta}{$\boldsymbol \beta$}
	\psfrag{beta(e)}{$X_{t,\tau}(M_l)$}
	\psfrag{Ext(beta,e)}[t][b]{$Ext_{t,\tau}(\boldsymbol \beta,M_l)$}
	\psfrag{e}[t][b]{$M_l$}
	\includegraphics[scale=0.55]{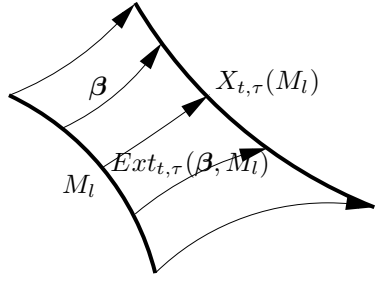}
    \end{center}
    \caption{Extrusion of line segment $M_l$ with respect to velocity field
    $\boldsymbol \beta$.}
    \label{fig:extrusion}
\end{figure}

$0$-forms can be identified with real valued functions. Recall the definition of the
directional derivative for a smooth function $f: \Omega \mapsto \mathbb R$:
\begin{eqnarray*}
    (\boldsymbol{\beta} \cdot
    \grad
    f)(x,t):=\lim_{\tau \to t} \frac{f(X(x,t,\tau))-f(x)}{\tau-t}\,.
\end{eqnarray*}
The Lie derivative $L_{\boldsymbol{\beta}}$ of $l$-forms $\omega \in
\mathcal{F}^l(\Omega)$ generalizes this. Note that the point evaluation of $0$-forms
is replaced with evaluation on $l$-dimensional oriented sub-manifolds $M_l \in
\mathcal S^l(\Omega)$. Thus the Lie derivative of a $l$-form $\omega$ is
\begin{eqnarray}
  <L_{\boldsymbol \beta} \omega(t),M_l>:=\lim_{ \tau \to t}
  \frac{<\omega,X_{t,\tau}(M_l)> - <\omega,M_l>}{\tau-t}\,.
\end{eqnarray}
In terms of the pullback $X^*_{t,\tau}$ defined by
\begin{eqnarray}
  < X^*_{t,\tau}\omega,M_l>:=<\omega,X_{t,\tau}(M_l)>
\end{eqnarray}
we can also write 
\begin{equation}
  L_{\boldsymbol \beta} \omega(t):=\lim_{ \tau \to t}
  \frac{X^*_{t,\tau}\omega - \omega}{\tau-t}\,.
  \label{eq:LieDefinition}
\end{equation}
Following \cite{BOS04a} the extrusion $Ext_{t,\tau}(\boldsymbol \beta,M_l)= \{
X_{t,s}(x) \,: \, t\leq s \leq \tau ,\, x \in M_l\} $ is the union of flux lines
emerging at $M_l$ running from $t$ to $\tau$ (Figure \ref{fig:extrusion}). The
orientation of $M_l$ induces an orientation of the extrusion
$Ext_{t,\tau}(\boldsymbol \beta,M_l)$ such that the boundary is:
\begin{equation} \partial Ext_{t,\tau}(\boldsymbol \beta, M_l)=
    X_{t,\tau}(M_l) - M_l -Ext_{t,\tau}( {\boldsymbol \beta}, \partial M_l)\,.
    \label{eq:ExtrusionBnd}
\end{equation} 
Plugging this into the definition of the Lie derivative
(\ref{eq:LieDefinition}) we get by means of Stokes' theorem
\begin{equation}
    \begin{aligned}
	<L_{\boldsymbol \beta} \omega(t),M_l> &= \lim_{\tau \to t} \frac{<\omega 
	,\partial Ext_{t,\tau}(\boldsymbol \beta, M_l)> + <\omega,Ext_{t,\tau}(\boldsymbol \beta,
	\partial M_l)>}{\tau-t} \\
	&{=} \lim_{\tau \to t}\frac{<\text{d} \omega,Ext_{t,\tau}(\boldsymbol \beta,
	M_l)> + <\omega, Ext_{t,\tau}(\boldsymbol \beta, \partial
	M_l)>}{\tau-t}\,.
	\label{eq:LieExtrusion}
    \end{aligned}
\end{equation}
The contraction operator is defined as the limit of the dual of the extrusion:
\begin{eqnarray}
    <i_{\boldsymbol \beta}\omega(t),M_l>:=\lim_{\tau \to t} \frac{<\omega,
      Ext_{t,\tau}(\boldsymbol \beta, M_l)>}{\tau - t}\;,
\end{eqnarray}
and we recover from (\ref{eq:LieExtrusion}) Cartan's formula \cite[page 142, prop.
5.3]{LangBook} for the Lie derivative:
\begin{equation}
  \begin{aligned}
    \label{eq:Cartan}
    L_{\boldsymbol \beta}\omega = i_{\boldsymbol \beta} \text{d}
    \omega + \text{d} i_{\boldsymbol \beta}\omega\;.
  \end{aligned}
\end{equation}
For $0$-forms the second term vanishes, for $n$-forms the first one.

For time dependent  
$l$-forms $\omega(t) \in \mathcal F^l(\Omega)$ the limit value:
\begin{eqnarray}
    D_{\boldsymbol \beta} \omega(t):=\lim_{ \tau \to t}
    \frac{X^*_{t,\tau}\omega(\tau) - \omega(t)}{\tau-t} 
    \label{eq:MatDefinition}
\end{eqnarray}
corresponding to \eqref{eq:LieDefinition} is the rate of change of the action of $l$-forms in
moving media, hence a material derivative \cite[page 62]{GUR81}.
We deduce:
\begin{equation}
    \begin{aligned}
    D_{\boldsymbol \beta} \omega(t) &= \lim_{ \tau \to t}
    \frac{X^*_{t,\tau}\omega(\tau) - X^*_{t,\tau}\omega(t)}{\tau-t}  
    + \lim_{ \tau \to t} \frac{X^*_{t,\tau}\omega(t)-
    \omega(t)}{\tau-t} \\
    &= \frac{\partial }{\partial t} \omega(t)+L_{\boldsymbol
    \beta}\omega(t)\,.
    \label{eq:material_deriv}
    \end{aligned}
\end{equation}

\begin{remark}
    Since the exterior derivative and the Lie derivative commute we have:
    \begin{eqnarray}
	D_{\boldsymbol \beta} d = d D_{\boldsymbol \beta}\,. 
    \end{eqnarray}
    As a consequence closed forms remain closed when they are advected by the
    material derivative.   
    If $D_{\boldsymbol \beta} \omega(t) = 0$ and 
    $d \omega(0) = 0$ then this property is preserved:
    $d \omega(t)=0,\, \forall t$.
    \label{remark:closed}
\end{remark}

\begin{remark}
    For $1$-forms $\omega \in \mathcal F^1(\Omega)$ with $\Omega \subset \mathbb R^3$
    with vector proxy $\mathbf u=\operatorname{v.p.}(\omega)$ formulas
    \eqref{eq:material_deriv} and \eqref{eq:Cartan} give a general convective term
    \begin{eqnarray*}
	\label{eq:Lie1F}
	\operatorname{v.p.}(L_{\boldsymbol \beta} \omega)= 
	-\boldsymbol \beta \times \curl {\mathbf u} + 
	\nabla (\boldsymbol \beta \cdot{\mathbf u})\;,
    \end{eqnarray*}
    while for $2$-forms $\omega \in \mathcal F^2(\Omega)$ with vector proxy
    $\mathbf u=\operatorname{v.p.}(\omega)$ the convective term reads 
    \begin{eqnarray*}
	\label{eq:Lie2F}
	\operatorname{v.p.} (L_{\boldsymbol \beta} \omega) = 
	\boldsymbol \beta( \Div {\mathbf u})+\curl(
	{\mathbf u}\times \boldsymbol \beta ) \;.
    \end{eqnarray*}
    For closed $2$-forms $\omega \in \mathcal F^l(\Omega)$, $\Omega \subset \mathbb R^3$
    with vector proxy $\mathbf u=\operatorname{v.p.}\omega$, which satisfies $\Div \ubf =0$, we find
    \begin{eqnarray*}
	\operatorname{v.p.}(D_{\boldsymbol \beta} \omega)=
	\partial_t {\mathbf u}+ 
	\curl(  {\mathbf u}\times \boldsymbol \beta ) \;.
    \end{eqnarray*}
\end{remark}

\subsection{Whitney forms}
\label{sec:discr-diff-forms}

We rely on a version of discrete exterior calculus spawned by using discrete
differential forms. We equip $\Omega$ with a simplicial triangulation $\Omega_h$. The
sets of subsimplices of dimension $l$ are denoted by $S_h^l(\Omega_h)$ and have
cardinality $N_l:=|S_h^l|$. In words, $S_h^0(\Omega_h)$ is the set of vertices of
$\Omega_{h}$, $S_h^1(\Omega_h)$ the set of edges, etc. We write $\mathcal
W^l(\Omega_h)\subset \mathcal F^l(\Omega)$ for the space of Whitney $l$-forms defined
on $\Omega_h$, see \cite[Sect.~23]{BOS05c} and \cite[Sect.~3]{AFW06}. From the finite
element point of view a Whitney $l$-form $\omega_h \in \mathcal W^l(\Omega_h)$ is
a linear combination
\begin{eqnarray}
    \omega_h=\sum_{i=1}^{N_l} b^l_i(x) \omega_i 
\end{eqnarray}
of certain basis functions associated with subsimplices of $\Omega_{h}$. The basis
functions $b^l_i$ have compact support and vanish on all cells that have trivial
intersection with subsimplex $s^l_i\in S_h^l(\Omega_h)$.  Locally, e.g. in each cell
$s^n_k$ with $\bar{s}_i^l\cap \bar{s}_n^k=\bar{s}_i^l$, the basis functions $b_i^l$
can be expressed in terms of the barycentric coordinate functions $\lambda_{j}$ and
whose gradients $d \lambda_i$. The barycentric coordinate functions $\lambda_j$ are
associated with vertices $s_j^0$ of cell $s_k^n$, e.g. $\bar{s}_j^0 \cap
\bar{s}^n_k
=\bar{s}_j^0$. If $I=(I_0,\dots I_l)$ denotes the index set of vertices $s_{I_0}^0,
\dots s_{I_l}^0$ of subsimplex $s_i^l$ then
\begin{eqnarray}
    b^l_{i|_{s^n_k}}= \sum_{j=0}^l (-1)^j \lambda_{I_j} \bigwedge_{k=1, k\neq j}^l
    d\lambda_{I_j}\,.    
\end{eqnarray}
The construction of these spaces ensures that the traces $\imath^{\ast} 
\omega_h$ on element boundaries
are continuous.  
The evaluation maps $(l^l_i)_{i=1}^{N_l}: \mathcal F^l \mapsto \mathbb
R^{N_l}$, $l_i^l(\omega):=<\omega,s^l_i>, s_i^l \in \mathcal S_h^l$
are the corresponding degrees of freedoms, which give rise to the usual nodal
interpolation operators $\Pi^l: \mathcal F^{l}(\Omega) \mapsto \mathcal
W^l(\Omega_h)$, also called DeRham maps \cite[Sect.~3.1]{HIP02},
\begin{eqnarray}
    \label{eq:nodalinterpol}
    \Pi^l \omega = \sum_{i=1}^{N_l}l_i^l(\omega) b^i_l.
\end{eqnarray}

\begin{remark}
  Discrete differential forms are available for more general triangulations
  \cite[Sect.~25]{BOS05c}. The Whitney forms feature piecewise linear vector
  proxies, but can be extended to schemes with higher polynomial degree
  \cite[Sect.~3.4]{HIP02}. For the sake of simplicity we are not going to 
  discuss those. 
\end{remark}

\section{Direct and adjoint semi-Lagrangian Galerkin formulation}
\label{sec:deriv_sLG}

There are two basically different strategies to discretize material derivatives. We
either approximate the difference quotient (\ref{eq:MatDefinition}) directly or use
the formulation (\ref{eq:material_deriv}) in terms of partial time derivative of
forms and Lie derivative. The former policy is referred to as Lagrangian and in our
setting of a single fixed mesh it is known as \emph{semi-Lagrangian}. It is explored
in this section. The latter is called \emph{Eulerian} and reduces to a method of lines
approach based on a spatial discretization of the Lie derivatives. Its investigation
is postponed to Sect.~\ref{sec:euler-appr-discr}.

We begin the study of the semi-Lagrangian method with a look at the pure transport
problem: given some $\phi(t) \in \mathcal{F}^{l}(\Omega)$ and $\omega(0) \in \mathcal
W^l(\Omega_h)$ find $\omega(t) \in \mathcal{F}^{l}(\Omega)$, such that
\begin{eqnarray}
    \label{eq:transport1}
    D_{\boldsymbol \beta}\omega=\phi\,.
\end{eqnarray}

Here and below we assume in the following that on the boundary of $\Omega$ the normal
components of the vector field ${\boldsymbol \beta}$ vanish, hence
$X(\Omega,\cdot,\cdot)=\Omega$. This is a restrictive but very common assumption for
semi-Lagrangian methods. Only a few semi-Lagrangian methods like ELLAM \cite{ELLAM}
can handle a non-vanishing velocity field on the boundary.
 
There are two different variational formulations for the transport problem
\eqref{eq:transport1}. The {\sl direct} variational formulation reads: given $\phi(t)
\in \mathcal{F}^{l}(\Omega)$ and initial data $\omega(0)\in \mathcal F^l(\Omega)$
find $\omega(t) \in \mathcal{F}^l(\Omega)$, such that
\begin{eqnarray}
    \int_{\Omega} \lim_{ \tau \to t}
    \frac{X^*_{t,\tau}\omega(\tau) - \omega(t)}{\tau-t} \wedge
    \eta= \int_{\Omega}\phi(t)
    \wedge \eta\,, \quad \forall \eta \in \mathcal{F}^{n-l}(\Omega)\,.
    \label{def:direct_var}
\end{eqnarray}
The {\sl adjoint} variational formulation for \eqref{eq:transport1} reads: given
$\phi(t) \in \mathcal{F}^{l}(\Omega)$ and initial data $\omega(0)\in \mathcal
F^l(\Omega)$ find $\omega(t) \in \mathcal{F}^l(\Omega)$, such that
\begin{eqnarray}
  \int_{\Omega} \lim_{ \tau \to t}
  \frac{ \omega(\tau) \wedge X^*_{\tau,t}\eta- \omega(t)\wedge
    \eta }{\Delta t} = \int_{\Omega}\phi(t)
  \wedge \eta\,, \quad \forall \eta \in \mathcal{F}^{n-l}(\Omega)\,.
  \label{def:weak_var}
\end{eqnarray}

\begin{lemma}
  \label{lem:31}
    The direct (\ref{def:direct_var}) and the adjoint (\ref{def:weak_var})
    variational formulations are equivalent.
\end{lemma}
\begin{proof}
    Let $\omega \in \mathcal F^l(\Omega)$ and $\eta \in \mathcal
    F^{n-1}(\Omega)$. Since $X_{t,s}(\Omega)=\Omega$ and $X_{\tau,t} \circ
    X_{t,\tau}=id$:
    \begin{eqnarray*}
	\int_{\Omega} X^*_{t,\tau} \omega(\tau) \wedge \eta  
	 =
	\int_{X_{\tau,t}(\Omega)} X^*_{t,\tau}\omega(\tau) \wedge \eta 
	 = \int_{\Omega} 
	\omega \wedge X^*_{\tau,t} \eta \;,
    \end{eqnarray*}
    and the equivalence follows.
\end{proof}

\begin{corollary}
    \label{cor:dualLie}
    Under the assumption $X_{t,\tau}(\Omega)=\Omega, \forall t,\tau$: 
    \begin{eqnarray*}
	\int_{\Omega}L_{\boldsymbol \beta} \omega\wedge \eta
	= - \int_{\Omega} \omega \wedge L_{\boldsymbol
	\beta}\eta\,, \quad \omega \in \mathcal F^l(\Omega),\eta \in \mathcal
	F^{n-l}(\Omega)\,.
    \end{eqnarray*}
\end{corollary}
Replacing the limits in (\ref{def:direct_var}) and (\ref{def:weak_var})
with a finite difference quotient yields semi-discrete timestepping schemes:
Given $\phi(t) \in \mathcal{F}^{l}(\Omega)$ and $\omega(t-\Delta t) \in
\mathcal{F}^l(\Omega)$ find $\omega(t) \in \mathcal{F}^l(\Omega)$, such that
\begin{eqnarray}
    \int_{\Omega} 
    \frac{\omega(t)- X^*_{t,t-\Delta t}\omega(t-\Delta t)}{\Delta t} \wedge
    \eta = \int_{\Omega}\phi(t)
    \wedge \eta\,, \quad \forall \eta \in \mathcal{F}^{n-l}(\Omega)
    \label{def:semi_direct_var}
\end{eqnarray}
and given $\phi(t) \in \mathcal{F}^{l}(\Omega)$ and $\omega(t-\Delta t) \in
\mathcal{F}^l(\Omega)$ find $\omega(t) \in \mathcal{F}^l(\Omega)$, such that
\begin{eqnarray}
    \int_{\Omega}
    \frac{\omega(t)\wedge
    \eta- \omega(t-\Delta t) \wedge X^*_{t-\Delta t,t}\eta}{\Delta t}  = \int_{\Omega}\phi(t)
    \wedge \eta\,, \quad \forall \eta \in \mathcal{F}^{n-l}(\Omega)\,.
    \label{def:semi_weak_var}
\end{eqnarray}
Restricting the semi-discrete
formulation to the discrete spaces $\mathcal W^l(\Omega_h)$, we end up with the
following direct and adjoint schemes:
Given $\phi(t) \in \mathcal{F}^{l}(\Omega)$ and $\omega(t-\Delta t) \in
\mathcal{W}^l(\Omega_{h})$ find $\omega_h(t) \in \mathcal{W}^l(\Omega_h)$, such that
\begin{eqnarray}
    \int_{\Omega} 
    \ast \frac{\omega_h(t)-X^*_{t,t-\Delta t}\omega_h(t-\Delta t)}{\Delta t} \wedge
    \tilde{\eta}_h =  \int_{\Omega} \ast\phi(t)
    \wedge \tilde{\eta}_h\,, \quad \forall \tilde{\eta}_h \in
    \mathcal{W}^{l}(\Omega_h)
    \label{def:semi2_direct_var}
\end{eqnarray}
and find $\omega_h(t) \in \mathcal{W}^l(\Omega_h)$, such that
\begin{multline}
    \int_{\Omega}
    \frac{\ast \omega_h(t)\wedge
    \tilde{\eta}_h- \ast \omega_h(t-\Delta t) \wedge
    X^*_{t-\Delta t,t} \tilde{\eta}_h}{\Delta t} \\ = \int_{\Omega} \ast \phi(t)
    \wedge \tilde{\eta}_h\,, \quad \forall \tilde{\eta}_h \in
    \mathcal{W}^{l}(\Omega_h)\,.
    \label{def:semi2_weak_var}
\end{multline}

A significant advantage of the semi-Lagrangian approach is the preservation of
closeness (see Remark \ref{remark:closed}). This property is important in many
physical applications. The adjoint semi-Lagrangian timestepping fulfils this in a
weak sense.
\begin{definition}
  A form $\omega \in \mathcal F^l(\Omega)$ is weakly closed if $\int_\Omega
  \ast \omega
  \wedge d \psi = 0, \forall \psi \in \mathcal F^{l-1}(\Omega)$.
\end{definition}
For $\phi=0$ the adjoint semi-Lagrangian timestepping boils down to
\begin{eqnarray*}
  \int_{\Omega} \ast \omega_h(t)\wedge d \psi_h = 
  \int_{\Omega} \ast \omega_h(t-\Delta t)\wedge X^*_{t-\Delta t,t} d \psi_h
  \quad \forall \psi_h \in
  \mathcal{W}^{l-1}(\Omega_h)\,.
\end{eqnarray*}
As the exterior derivative and pullback commute we conclude
\begin{eqnarray*}
  \int_{\Omega} \ast \omega_h(t)\wedge d \psi_h = 
  \int_{\Omega} \ast \omega_h(t-\Delta t)\wedge d \widetilde{\psi_h}\,,
\end{eqnarray*}
with $\widetilde{\psi_h}= X^*_{t-\Delta t,t}\psi_h$.
Hence $\omega_h(t)$ is weakly closed if $\omega_h(t-\Delta t)$ is
weakly closed.
\label{remark:weakclosed}

\begin{remark}
  \label{remark:stab_const_cont}
  Another advantage of semi-Lagrangian methods is the straightforward stability
  analysis for the homogeneous transport problem $D_{\boldsymbol \beta}\omega_h=0$.
  The corresponding direct scheme reads: given $\omega_h(\tau) \in \mathcal
  W^l(\Omega_h)$, find $\omega_h(t) \in \mathcal W^l(\Omega_h)$
    \begin{eqnarray*}
	\int_\Omega \ast \omega_h(t) \wedge \tilde \eta_h = \int_\Omega
	\ast X^{*}_{t,t-\Delta t} \omega_h(t-\Delta t) \wedge \tilde \eta_h\,, \quad \forall
	\tilde \eta_h
	\in \mathcal W^l(\Omega_h)\,.
    \end{eqnarray*}
    Testing with $\omega_h$ we immediately get the stability estimate:
    \begin{eqnarray*}
      \| \omega_h(t) \|_{0}  
      &\leq& \| X^*_{t,t-\Delta t}\omega_h(t-\Delta t) \|_{0}\\
      &\leq& C \| \omega_h(t-\Delta t) \|_{0}
    \end{eqnarray*}
    with $C=C(\boldsymbol \beta, \Delta t)>0$. For $\Omega \subset \mathbb R^3$ and Euclidean
Hodge the concrete constants can be deduced from vector proxy representations
of the pullback \cite[p. 245]{HIP02}. The stability estimates then read:
\begin{equation*}
    \begin{aligned}
	\| \omega_h(t) \|_{0}&\leq \sup_{x \in
	\Omega}\,|\det(DX_{t-\Delta t,t}(x))|^{\frac12}\,\|\omega_h(t-\Delta t)\|_{0}, \quad
	\omega_h \in \mathcal W^0(\Omega_h) ; \\ 
	\| \omega_h(t) \|_{0}&\leq \sup_{x \in
	\Omega}\, \left(|(\det(DX_{t-\Delta t,t}(x))|^{\frac12}\,
	\|DX_{t,t-\Delta t}(x)\|\right) \,\|\omega_h(t-\Delta t)\|_{0}, \quad
	\omega_h \in \mathcal W^1(\Omega_h); \\ 
	\| \omega_h(t) \|_{0}&\leq \sup_{x \in
	\Omega}\,\left(|\det(DX_{t,t-\Delta t}(x))|^{\frac12}\, \|DX_{t-\Delta
	t,t}(x)
	\|\right)\,\|\omega_h(t-\Delta t)\|_{0}, \quad \omega_h \in \mathcal
	W^2(\Omega_h); \\ 
	\| \omega_h(t) \|_{0}&\leq \sup_{x \in
	\Omega}\,|\det(DX_{t,t-\Delta t}(x))|^{\frac12}\,\|\omega_h(t-\Delta t)\|_{0}, \quad
	\omega_h \in \mathcal W^3(\Omega_h), 
    \end{aligned}
\end{equation*}
where $\|\cdot\|$ is the Euclidean matrix norm. Stability will depend in
general on the distortion effected by $X_{t,t-\Delta t}$. For $l=1,2$ even the
usual assumption $\Div \boldsymbol \beta=0$ does not guarantee stability. 
\end{remark}

It is important to note that both the stability estimates and the preservation
of weak closedness as previously stated assume exact evaluation of the
bilinear forms $\int_\Omega
\ast X^{*}_{t,t-\Delta t} \omega_h(x,t-\Delta t) \wedge \tilde \eta_h(x) dx$ and $\int_\Omega
\ast \omega_h(x,t-\Delta t) \wedge X^{*}_{t-\Delta t,t} \tilde \eta_h(x) dx$.
Any non-trivial problem will require the use of additional approximations. We
propose the following three approximation steps. 
\begin{enumerate}
\renewcommand{\labelenumi}{(\roman{enumi}) }
\item First we use some numerical integrator to determine approximations $\tilde
  X_{t,t'}(s^0_i)$ to the exact evolutions $X_{t,t'}(s^0_i)$ of vertices $s^0_i$.
  Note that the direct scheme requires to solve \eqref{eqn:charac_ode} backward in
  time, while for the adjoint scheme we calculate forward in time. The simplest
  numerical integrator for this first step is the forward Euler scheme, which gives
  $\tilde X_{t,t'}(s^0_i)=s_i^0+(t'-t) \boldsymbol \beta(s_i^0)$.
\item Then we approximate the
  evolution $X_{t,t'}(x)$ of arbitrary points $x$ as the convex combination
  of the approximated evolutions of surrounding vertices:
  \begin{eqnarray}
    X_{t,t'}(x)\approx \bar{X}_{t,t'}(x):=\sum_{i=1}^{N_0} \tilde{X}_{t,t'}(s^0_i)
    \lambda_i(x)\,, 
  \end{eqnarray}
  This yields a piecewise linear, continuous approximation $\bar{X}_{t,t'}(x)$
  of the exact evolution $X_{t,t'}(x)$.

  Note that for small time steps and Lipschitz continuous $\boldsymbol \beta$ the
  approximating flow $\bar X_{t,t'}$ maps the simplicial triangulation $\Omega_h$
  onto a new simplicial triangulation $ \bar \Omega_h$, such that each subsimplex
  $s_i^l$ of $\Omega_h$ has an image $\bar s_i^l=: \bar X_{t,t'}(s_i^l)$ in $\bar
  \Omega_h$.
\item Finally, we apply the nodal interpolation operators of discrete forms
  \cite[Sect.~3.3]{HIP02} to map the transported discrete forms $\bar{X}^*_{t,t-\Delta
    t}\omega_h(t-\Delta t)$ and $\bar{X}^*_{t-\Delta t,t}\tilde\eta_h$ onto the space
  of discrete forms.
\end{enumerate}

\begin{figure}[h]
    \begin{center}
	\psfrag{a_i}{$s^0_i=X_{t,t}(s^0_i)$}
	\psfrag{a_k}{$s^0_k$}
	\psfrag{a_l}{$s^0_l$}
	\psfrag{a_m}{$s^0_m$}
	\psfrag{T_1}{$s^n_1$}
	\psfrag{T_2}{$s^n_2$}
	\psfrag{T_3}{$s^n_3$}
	\psfrag{T_4}{$s^n_4$}
	\psfrag{X(a_i)}{$X_{t,t-\Delta t}(s^0_i)$}
	\includegraphics[width=0.7\linewidth]{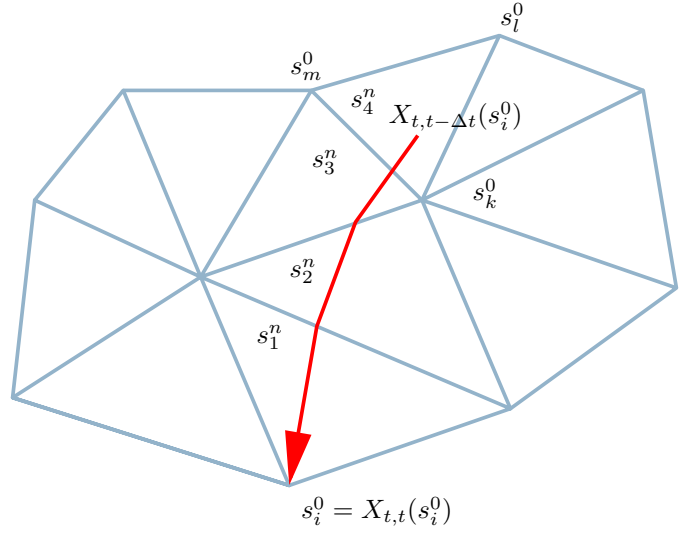}
    \end{center}
    \caption{To determine the location of $\bar{X}_{t,t-\Delta t}(s^0_i)$, we move
    backward along the
    trajectory $\bar{X}_{t,\cdot}(s^0_i)$ starting from $s^0_i$ and identify the crossed
    elements $s^0_1,s^0_2,s^0_3$ and $s^0_4$. In this case 
    $p^0_{ik},p^0_{il}$ and $p^0_{im}$ are the only non-zero entries in
    the $i$--th row of ${\mathbf P}_{t,t-\Delta t}^0$.}
    \label{fig:transport_vertex}
\end{figure}

Let us illustrate these three steps for the case of $0$-forms and $1$-forms: To
determine now the interpolant of a discrete transported $0$-form $\bar
X^{*}_{t,t'}\omega_h, \omega_h \in \mathcal W^0(\Omega_h)$, by linearity it is enough
to consider the basis functions. Since $<\Pi^0\omega,s^0_i>=<\omega,s^0_i>$ for all
vertices $s^0_i$ and $0$--forms $\omega$ the matrix operator ${\mathbf P}_{t,t'}^0$
with entries
\begin{eqnarray}
    p^0_{ij}:=l_i^0(\bar{X}_{t,t'})=\lambda_j(\bar{X}_{t,t'}(s^0_i))=\lambda_j(\tilde{X}_{t,t'}(s^0_i))
    \label{def:P0_entries}
\end{eqnarray}
maps the expansions coefficients of $\omega_h$ to the coefficients of
$\Pi^0\bar{X}^*_{t,t'}\omega_h$. This means that in each time step we not only
need to determine the points $\tilde{X}_{t,t'}(s^0_i)$ but also the location within the
triangulation $\Omega_h$. To
find the element, in which $\bar{X}_{t,t'}(s^0_i)$ is located, we trace the
path of the
trajectory from one element to the next. Based on this data the matrix entries
(\ref{def:P0_entries}) can be assembled element by element (see fig.
\ref{fig:transport_vertex} for the direct scheme).

\begin{figure}[h]
    \begin{center}
	\psfrag{a_i}{$s^0_j$}
	\psfrag{a_j}{$s^0_k$}
	\psfrag{a_l}{}
	\psfrag{T_1}{$s^n_1$}
	\psfrag{T_2}{$s^n_2$}
	\psfrag{T_3}{$s^n_3$}
	\psfrag{X(a_i)}{$X_{t,t-\Delta t}(s^0_j)$}
	\psfrag{X(a_j)}{$X_{t,t-\Delta t}(s^0_k)$}
	\includegraphics[width=0.7\linewidth]{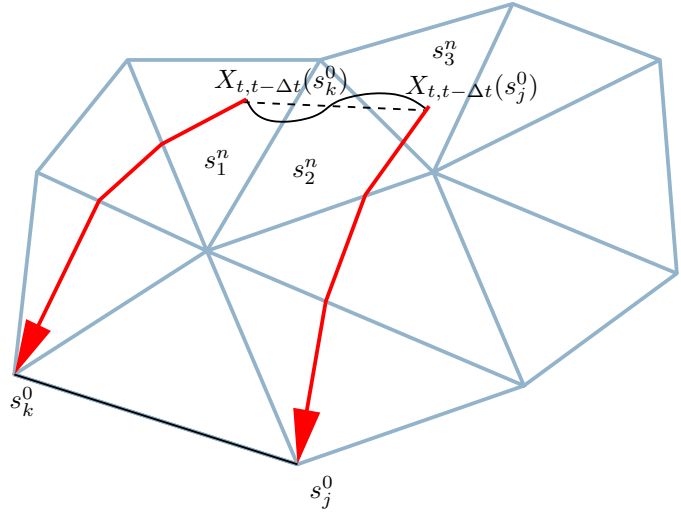}
    \end{center}
    \caption{The transported edge ${X}_{t,t-\Delta t} s_i^1$ (black curved line) 
    is approximated by a straight line $\bar {X}_{t,t-\Delta t} s_i^1$(black
    dashed line). In the case
    depicted here all basis functions associated with edges $s_{i'}^1$ of
    elements $s^n_1$,
    $s^n_2$ and $s^n_3$ yield a non-zero entry $p^1_{i, i'}$. }
    \label{fig:transport_edge}
\end{figure}

\begin{figure}[h]
    \begin{center}
	\psfrag{a}{$a$}
	\psfrag{b}{$b$}
	\psfrag{e1}{$s_{j'}^0$}
	\psfrag{T}{$s^n_{\cdot}$}
	\psfrag{e2}{$s_{k'}^0$}
	\psfrag{e}{$s_{i'}^1$}
	\includegraphics[width=0.7\linewidth]{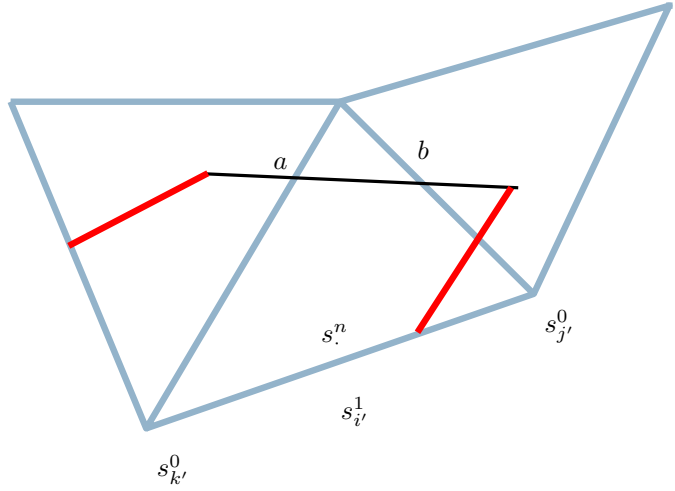}
    \end{center}
    \caption{The line $[a,b]$ is the intersection of the approximation of the
    transported edge with element $s_{\cdot}^n$.}
    \label{fig:sketch_edge}
\end{figure}

The advantage of the second approximation step, the linear interpolation, is
elucidated by the treatment
of $1$-forms. The interpolation of a transported discrete $1$--form $\bar{X}^*_{t,t'}\omega_h$
is determined by the interpolation of transported basis forms and
the condition
$<\Pi^1\bar{X}^*_{t,t'}\omega_h-\bar{X}^*_{t,t'}\omega_h,s_i^1>=0$ for all edges
$s_j^1 \in \mathcal S_h^1$. This defines a matrix ${\mathbf P}_{t,t'}^1$, mapping
the expansion coefficients of $\omega_h$ to those of
$\bar{X}^*_{t,t'}\omega_h$. The matrix entries 
\begin{eqnarray}
    \label{eq:Pmat1F}
    p_{ i i'}^1:=l_i^1(\bar{X}^*_{t,t'}b_{i'}^1)=<\bar{X}^*_{t,t'}b_{i'}^1,s_i^1>\,,
\end{eqnarray}
are line integrals along the straight line from $\tilde{X}_{t,\tau}(s^0_j)$ to
$\tilde{X}_{t,\tau}(s^0_{k})$, if $s_i^1$ is the edge connecting the vertices $s^0_j$
and $s^0_k$ (see fig. \ref{fig:transport_edge} for the direct scheme).    
To determine the
entries of the $i$-th row, we trace the straight line from
$\bar{X}_{t,t'}(s_j^0)$ to
$\bar{X}_{t,t'}(s_k^0)$ and calculate for each crossed element the line
integrals for the attached basis functions. If e.g. the line crosses an
element $s_{\cdot}^n$ with edge $s_{i'}^1$ from point $a$ to point $b$ (see fig.
\ref{fig:sketch_edge} for the direct scheme), then the
element contribution to $p^1_{ i i'}$ is:
\begin{eqnarray*}
    \int_{\bar{X}_{t,t'}(s_i^1 \cap s_{\cdot}^n)}b_{i'}^1 =
    \int_{[a,b]}b_{i'}^1
    =
    \lambda_{j'}(a)\lambda_{k'}(b)-
    \lambda_{k'}(a)\lambda_{j'}(b)\,,
\end{eqnarray*}
where $s_{j'}^0$ and $s_{k'}^0$ are the terminal vertices of edge $s_{i'}^1$.

Observe that we do not need to calculate additional approximate solutions of
\eqref{eqn:charac_ode} to determine the matrix entries \eqref{eq:Pmat1F}.  Piecewise
linear interpolation of the approximate evolution of vertices automatically gives
approximate evolutions of all subsimplices.

Although the approximate flow $\bar X_{t,t'}$ is not smooth, we can show that the
fully discrete adjoint scheme preserve the weak closedness property.  Is is easily
established that exterior derivative and interpolation of approximate pullbacks
commute \cite[Sect.~3.2]{HIP02}.
\begin{lemma}
    \label{lemma:commute}
    Let $\Pi^l: \mathcal F^l(\Omega) \mapsto \mathcal W^l(\Omega_h)$ be the standard
    nodal interpolation operator and $\bar{X}_{t,\tau}(x):=\sum_i \tilde X_{t,\tau}(x_i)
    \lambda_i(x)$, where the coefficients $\tilde X_{t,\tau}$ are approximate
    solutions of the characteristic ODE \eqref{eqn:charac_ode}, then we have:
    \begin{eqnarray*}
	d \Pi^{l} \bar{X}^*_{t,\tau}\omega = \Pi^{l+1} \bar{X}^*_{t,\tau} d
	\omega\,, \quad \forall \omega \in \mathcal W^{l}(\Omega_h)\,.
    \end{eqnarray*}
\end{lemma}

Finally we state the fully discrete Semi-Lagrangian timestepping schemes for
the general convection-diffusion problem \eqref{eq:cddf}:

\begin{definition}
  Given a temporal mesh $0=t_{0}<t_{1}<t_{2}<\dots$, the discrete semi-Lagrangian
  timestepping scheme for the direct variational formulation of convection-diffusion
  of differential forms is: Given $\omega_h(t_0) \in \mathcal W^l(\Omega_h)$ and
  $\phi \in \mathcal F^{l}(\Omega)$, find $\omega_h(t_k), k>0$ such that:
\begin{multline}
    \label{def:SL_cd_dir}
    \int_{\Omega} \ast \frac {\omega_h(t_{k+1})-
    \Pi^l \bar{X}_{t_{k+1},t_k}^ *\omega_h(t_k)}{t_{k+1}-t_k} \wedge
    \tilde\eta_h
    \\ + \int_{\Omega}
    \ast d \omega_h(t_{k+1})
    \wedge d \tilde \eta_h = \int_{\Omega} \ast \phi \wedge \tilde\eta_h\,, \quad
    \forall \tilde\eta_h \in
    \mathcal W^l(\Omega_h)\,.
\end{multline}
The discrete Semi-Lagrangian timestepping scheme for
the adjoint variational formulation of convection-diffusion of differential
forms is: Given $\omega_h(t_0) \in
\mathcal W^l(\Omega_h)$ and $\phi \in \mathcal F^l(\Omega)$, find
$\omega_h(t_k), k>0$ such that:
\begin{multline}
    \label{def:SL_cd_ad}
    \int_{\Omega} \frac {\ast \omega_h(t_{k+1})\wedge \tilde\eta_h
    - \ast \omega_h(t_k)\wedge \Pi^l
    \bar{X}_{t_k,t_{k+1}}^*\tilde\eta_h
    }{t_{k+1}-t_{k}}    \\ + \int_{\Omega}
    \ast d \omega_h(t_{k+1})
    \wedge d \tilde\eta_h  =
    \int_{\Omega}\ast \phi \wedge \tilde\eta_h\,, \quad \forall \tilde\eta_h \in \mathcal
    W^l(\Omega_h)\,.
\end{multline}
\end{definition}

\begin{remark}
  The proposed sequence of approximation steps allows for a comprehensive
  stability and convergence analysis of fully discrete schemes. Stability of the
  fully discrete schemes in $\Omega=\mathbb R^3$ for example follows from  
  remark \ref{remark:stab_const_cont}. We only need to control perturbations due to some
  numerical integrator, the linear interpolation and interpolation of forms.
  We will address this in some future work.
\end{remark}

\begin{remark}
  The standard finite element approach for approximating bilinear forms like $
  \int_\Omega \ast {X}^{*}_{t,t'} \omega_h \wedge \tilde \eta_h$ would be quadrature.
  But since the product of ${X}^{*}_{t,t'} \omega_h$ and $\tilde\eta_h$ is smooth
  only on intersections $s^n_i\cap X_{t',t}(s^n_j)$, this entails finding all
  intersections of an element $s^n_i$ of the triangulation with any other element
  $X_{t',t}(s^n_j)$ of the transported triangulation, and apply some quadrature
  there. We conclude that such an approach is even more complex than methods based in
  interpolation. Further, it is totally unclear how such an approximation affects
  stability and conservation of closedness.
\end{remark}

\begin{remark}
    If we use in the scalar case $l=0$ a vertex based quadrature rule to
    approximate the element integrals $\int_{s_i^n} \Pi^0 \bar X^*_{t,\tau}
    \omega_h(\tau) \wedge \eta$, the resulting scheme corresponds to a
    Lagrange-Galerkin scheme from \cite{MPS88} using vertex based quadrature.
\end{remark}

\begin{remark}
  Notice that in both \eqref{def:SL_cd_dir} and \eqref{def:SL_cd_ad} the unknown
  $\omega_{h}(t_{k+1})$ is obtained by solving a standard discrete diffusion problem
  arising from a symmetric positive definite bilinear form, which are amenable to
  fast iterative solvers \cite{HIX06,AFW99,HIP99}.  In the words of \cite{JXU08}, the
  semi-Lagrangian method is ``solver-friendly''.
\end{remark}

\section{Eulerian approach: discrete Lie-derivatives}
\label{sec:euler-appr-discr}

The Eulerian approach treats the discretization in space and the
discretization in time separately. We therefore consider
the stationary version of the convection-diffusion equation
\eqref{eq:cddf}: 
\begin{align}
    \label{eq:cddf_stat}
    \begin{array}[c]{rcll}
    - \varepsilon (-1)^{l} 
    d\ast d \omega +  \ast
    L_{\boldsymbol \beta} \omega &=& \varphi \quad &
    \text{in }\Omega\subset\mathbb{R}^{n}\;,\\
    \imath^{\ast}\omega &=& 0 & \text{on }\partial\Omega\;.
  \end{array}
\end{align}
Since the approximation of the diffusion term is well-known \cite{HIP99c}, we
will discuss only the approximation of the convective terms in \eqref{eq:cddf_stat}.
We will pursue two different ways to obtain spatial discretization for the
Lie-derivative terms.  

\subsection{Standard scheme}
The standard Galerkin technique to approximate \\\mbox{$\int_{\Omega}\ast
  L_{\boldsymbol \beta}\omega_h \wedge \eta_h, \,\,\omega_h,\eta_h \in \mathcal
  W^l(\Omega_h)$} relies on elementwise quadrature. While evaluation of the part
\mbox{$\int_{\Omega}\ast i_{\boldsymbol \beta} d \omega_h \wedge \eta_h$} is standard
in the finite element context, we encounter some difficulties for the term
\mbox{$\int_{\Omega}\ast d i_{\boldsymbol \beta} \omega_h \wedge \eta_h$}, which
occurs for $l>0$. The term $d i_{\boldsymbol \beta}$ implicitly requires spatial
derivatives that fail to be square integrable for discrete differential forms. For
instance, for $1$-forms $\omega \in \mathcal F^l(\Omega)$ with smooth vector proxy
$\ubf:=\operatorname{v.p.}(\omega)$ we have:
\begin{align}
    \grad(\boldsymbol \beta \cdot \ubf)&= (D\boldsymbol \beta)^T \ubf +
    (D\ubf)^T \boldsymbol \beta \\
    &=  (D\boldsymbol \beta)^T \ubf- 
    (D\ubf)_{asym} \boldsymbol \beta+ (D\ubf)_{sym} \boldsymbol \beta\,, 
\end{align}
while for discrete $1$-forms $\omega_h \in \mathcal W^l(\Omega_h)$ with vector proxy
$\ubf_h $ on each element the symmetric part
$(D\ubf_h)_{sym}:=\frac12(D\ubf_h+(D\ubf_h)^T)$ vanishes. A scheme based on
elementwise quadrature will always miss $(D\ubf)_{sym} \boldsymbol \beta$.  Stated
differently, the exterior derivative $d i_{\boldsymbol \beta} \omega_h$ of the
contraction of a discrete form is not defined in a strong sense. But it is defined in
a weak sense: For any smooth $n-l$ form $\gamma \in \mathcal F^{n-l}(\Omega),
\imath^\ast \gamma=0$, and discrete differential form $\omega_h \in \mathcal
W^l(\Omega_h) $ elementwise integration by parts yields:
\begin{eqnarray}
    (-1)^{l} \int_\Omega i_{\boldsymbol \beta} \omega_h \wedge d \gamma =
    \sum_{i=1}^{N_n} \int_{s_i^n} d_hi_{\boldsymbol \beta} \omega_h \wedge \gamma -
    \sum_{i=1}^{N_n} \int_{\partial s_i^n} i_{\boldsymbol \beta}\omega_h \wedge
    \gamma\,,   
\end{eqnarray}
where $d_h$ is the restriction of $d$ to each element $s_i^n$. Since each face inside
$\Omega$ occurs twice in the right sum we can rewrite this term as sum of face
integrals over certain jump terms $[i_{\boldsymbol \beta}\omega_h]^{l-1}$. 
\begin{eqnarray}
    \label{eq:distr_def}
    (-1)^{l} \int_\Omega i_{\boldsymbol \beta} \omega_h \wedge d \gamma =
    \sum_{i=1}^{N_n} \int_{s_i^n} d_hi_{\boldsymbol \beta} \omega_h \wedge \gamma -
    \sum_{j=1}^{N_{n-1}} \int_{s_{j}^{n-1}} [i_{\boldsymbol
    \beta}\omega_h]^{l-1} \wedge
    \gamma\,,   
\end{eqnarray}
One is now tempted to use the right hand side of
\eqref{eq:distr_def} to define approximations for the convection terms. The
difficulty here are the face integrals, which are not well-defined for
$\gamma \in \mathcal W^{l}(\Omega_h)$. Motivated by finite volume
approximations, we replace in the case $\gamma \in \mathcal W^l(\Omega_h)$ in
the face integrals $\gamma$ with a flux function
$\hat{\gamma}=\hat{\gamma}(\gamma^l,\gamma^r)$ depending on the values
$\gamma^l$ and $\gamma^r$ on adjacent elements.

We use the characterization (\ref{eq:LieDefinition}) of the Lie derivative to determine such
consistent flux functions. 
\begin{definition}
  \label{def:bd}
    Let $\boldsymbol \beta$ be Lipschitz continuous, $\omega, \tilde\eta \in
    \mathcal F^l(\Omega)$ and $\Delta \tau>0$  then
    \begin{align}
	\label{eq:var_ufd}
	b_{\Delta \tau}(\omega ,\tilde\eta):= \int_\Omega
	\ast L_{\boldsymbol \beta}^{\Delta \tau} \omega\wedge \tilde\eta 
    \end{align}
    is the variational \emph{central} finite difference, with
    \begin{align*}
	L_{\boldsymbol \beta}^{\Delta \tau}
	\omega = \frac{X^{*}_{t,t+\Delta \tau} \omega-X^*_{t,t-\Delta
	\tau}\omega }{2 \Delta
	\tau}\,.
    \end{align*}
\end{definition}

\begin{figure}[h]
    \begin{center}
	\psfrag{ri}{$X_{t-\Delta \tau,t}(s^2_i)$}
	\psfrag{rj}{$X_{t-\Delta \tau,t}(s^2_j)$}
	\psfrag{bi}{$s^2_i$}
	\psfrag{bj}{$s^2_i$}
	\psfrag{a}{$a$}
	\psfrag{X(a)}{$X_{t-\Delta \tau,t}(a)$}
	\includegraphics[width=0.4\linewidth]{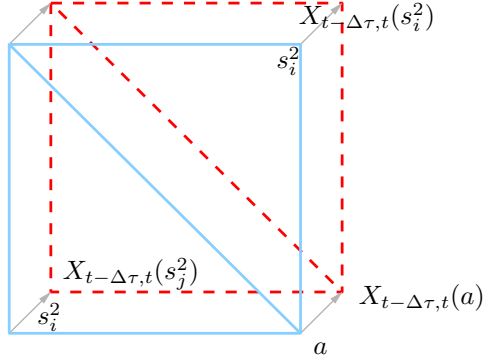}
    \end{center}
    \caption{Illustration for proof of Lemma \ref{lemma:standard1}. In the
    case of the upwind finite difference we have, that discrete $l$-forms
    $\Omega_h \in \mathcal W^l, l>0$ are discontinuous on edges of the blue
    triangulation $\Omega_h$. Their pullbacks $X^{*}_{t,t-\Delta \tau} \omega_h$ are
    discontinuous on edges of
    the red dashed triangulation $X_{t-\Delta\tau,t}(\Omega_h)$, the image of
    $\Omega_h$.
  }   
    \label{fig:limit_pic}
\end{figure}

We could have used one-sided finite differences as well.  Next we show that the limit
$\lim_{\Delta \tau \to 0} b_{\Delta \tau}(\omega_h, \tilde\eta_h), \omega_h,
\tilde\eta_h \in \mathcal W^l(\Omega_h)$ exists also for discrete differential forms.
\begin{lemma}
    \label{lemma:standard1}
    The limit $b_0(\omega_h,
	\tilde\eta_h):= \lim_{\Delta \tau \to 0} b_{\Delta \tau}(\omega_h, \tilde\eta_h),
    \omega_h,\tilde\eta_h \in \mathcal W^l(\Omega_h)$ exists, and we have the
    representation
    \begin{align}
	b_0(\omega_h,
	\tilde\eta_h)=\sum_{i=1}^{N_n} \int_{s_i^n} \ast \bar L_{\boldsymbol \beta}
	\omega_h \wedge \tilde\eta_h - \sum_{j=1}^{N^{\circ}_{n-1}}
	\int_{s_{J_j^{\circ}}^{n-1}}
	\frac12i_{\boldsymbol\beta}\left(\ast (\omega_h^+-
	\omega_h^-)\wedge (\tilde\eta_h^+ + \tilde\eta_h^-)\right),
	\label{eq:limit_standard}
    \end{align}
    with upwind and downwind traces $\omega_h^+$ and $\omega_h^-$ and the
    pointwise limit $\bar{L}_{\boldsymbol \beta}\omega_h(x)=\lim_{\Delta \tau \to
    0} L_{\boldsymbol \beta}^{\Delta \tau} \omega_h(x)$. $J^{\circ} \subset
    \{1,\dots N_{n-1}\}$ is the
    index set of faces, that are not on the boundary of $\Omega$. It's
    cardinality is $N^{\circ}_{n-1}$. Here we used
    \begin{equation}
	\int_{s_{k}^{n-1}} i_{\boldsymbol \beta}
	(\ast \omega_h^- \wedge \tilde\eta_h^+):=
	\lim_{\Delta
	\tau \to 0} \int_{Ext_{t-\Delta \tau,t}(\boldsymbol
	\beta,s_{k}^{n-1})}\ast
	\frac{X_{t,t-\Delta
	\tau}^* \omega_h}{\Delta \tau} \wedge \tilde \eta_h\,. 
    \end{equation}
\end{lemma}
\begin{proof}
    Let us first note that for all $x \not \in \partial s_{i}^n$ the pointwise
    limit $\bar{L}_{\boldsymbol \beta}\omega_h(x):=\lim_{\Delta \tau \to
    0} L_{\boldsymbol \beta}^{\Delta \tau} \omega_h(x)$ exists.

    We decompose the central difference into a sum of upwind and downwind finite
    differences:
    \begin{equation*}
	\begin{aligned}
	    b_{\Delta \tau}(\omega_h,\tilde\eta_h)=\frac12 \left(\int_{\Omega}
	    \ast \frac{X^{*}_{t,t+\Delta \tau}\omega_h - \omega_h}{\Delta \tau}
	    \wedge \tilde\eta_h + \int_{\Omega}
	    \ast \frac{\omega_h - X^{*}_{t,t-\Delta \tau}\omega_h}{\Delta \tau}
	    \wedge \tilde\eta_h\right) \;,
	\end{aligned}
    \end{equation*}
    and look first at the limit for the upwind finite difference. We split the
    integration over the domain $\Omega$ in a sum of integrals over patches
    $P_{i,j}(\Delta \tau):=s^n_{i}\cap
    X^{-1}_{t,t-\Delta \tau} (s^n_j)$ with smooth integrand:
    \begin{equation*}
	\begin{aligned}
	   \int_{\Omega}
	    \ast \frac{\omega_h - X^{*}_{t,t-\Delta \tau}\omega_h}{\Delta \tau}
	    \wedge \tilde\eta_h
	    & =\sum_{i,j=1}^{N_n} \int_{P_{i,j}(\Delta \tau)
	    } \ast \frac{\omega_h - X^{*}_{t,t-\Delta
	    \tau}\omega_h}{\Delta \tau}\wedge \tilde\eta_h :=\sum_{i,j=1}^{N_n} I_{i,j}(\Delta \tau) 
	\end{aligned}
    \end{equation*}
    and distinguish three cases, see Fig.~\ref{fig:limit_pic}:
    \begin{itemize}
	\item $i\neq j $ and $\bar s^n_j\cap \bar s^n_i = \emptyset$:
	    \begin{equation}
		\lim_{\Delta \tau \to 0} I_{i,j}(\Delta \tau)=0
	    \end{equation}
	\item $i\neq j $ and $\bar s^n_j\cap \bar s^n_i = s^{n-1}_k$: This
	    implies first that $s^{n-1}_k \not \subset \partial \Omega$.
	    If elements $s_i^n$ and $s_j^n$ share a face $s_k^{n-1}$ it is
	    clear that  
	    $P_{i,j}(\Delta \tau)=Ext_{t-\Delta \tau,t} (\boldsymbol
	    \beta, s_{k}^{n-1})+O(\Delta \tau^2)$ and we  have
	    \begin{equation}
		\begin{aligned}
		    \lim_{\Delta \tau \to 0} \int_{P_{i,j}(\Delta \tau)} \ast
		    \frac{\omega_h}{\Delta \tau} \wedge \tilde\eta_h& =\lim_{\Delta
		    \tau \to 0} \int_{Ext_{t-\Delta \tau,t}(\boldsymbol
		    \beta,s_{k}^{n-1})}\ast
		    \frac{\omega_h}{\Delta \tau} \wedge \tilde\eta_h \\
		    &= -\int_{s_{k}^{n-1}} i_{\boldsymbol \beta}
		    (\ast \omega_h^+ \wedge \tilde\eta_h^+)
		\end{aligned}
	    \end{equation}
	    and 
	    \begin{equation}
		\begin{aligned}
		    \lim_{\Delta \tau \to 0} \int_{P_{i,j}(\Delta \tau)}\ast \frac{X_{t,t-\Delta
		    \tau}^* \omega_h}{\Delta \tau} \wedge \tilde\eta_h&=\lim_{\Delta
		    \tau \to 0} \int_{Ext_{t-\Delta \tau,t}(\boldsymbol
		    \beta,s_{k}^{n-1})}\ast
		    \frac{X_{t,t-\Delta
		    \tau}^* \omega_h}{\Delta \tau} \wedge \tilde\eta_h \\
		    &= -\int_{s_{k}^{n-1}} i_{\boldsymbol \beta}
		    (\ast \omega_h^- \wedge \tilde\eta_h^+)
		\end{aligned}
	    \end{equation}
	    by the definition of contraction.
	\item $i=j:$
	    Since by assumption $\lim_{\Delta \tau \to 0} \ast
	    \frac{\omega_h - X_{t,t-\Delta
	    \tau}^* \omega_h}{\Delta \tau}$ is bounded on $P_{i,i}(\Delta
	    \tau)$ we get:
	    \begin{equation}
		\begin{aligned}
		    \lim_{\Delta \tau \to 0} \int_{P_{i,i}(\Delta \tau)}\ast
		    \frac{\omega_h - X_{t,t-\Delta
		    \tau}^* \omega_h}{\Delta \tau} \wedge \tilde\eta_h =& \\ 
		    \int_{s_i^n} \ast \bar L_{\boldsymbol \beta} \omega_h
		\wedge \tilde\eta_h\,. 
	    \end{aligned}
	    \end{equation}
    \end{itemize}
    After summation over all elements we obtain:
    \begin{equation}
	\begin{aligned}
	    \lim_{\Delta \tau \to 0} \int_{\Omega}
	    \ast \frac{\omega_h - X^{*}_{t,t-\Delta \tau}\omega_h}{\Delta \tau}
	    \wedge \tilde\eta_h = & \\
	    & \hspace{-3cm}\sum_{i=1}^{N_n} \int_{s_i^n} \ast \bar L_{\boldsymbol \beta}
	    \omega_h \wedge \tilde\eta_h - \sum_{j=1}^{N^{\circ}_{n-1}}
	    \int_{s_{J^{\circ}_j}^{n-1}}
	    i_{\boldsymbol\beta}\left( \ast (\omega_h^+-
	    \omega_h^-)\wedge \tilde\eta_h^+\right)
        \end{aligned}
    \end{equation}
    for the upwind finite difference.
    A similar argument for the downwind finite difference shows:
    \begin{equation}
	\begin{aligned}
	    \lim_{\Delta \tau \to 0} \int_{\Omega}
	    \ast \frac{X^{*}_{t,t+\Delta \tau}\omega_h-\omega_h}{\Delta \tau}
	    \wedge \tilde\eta_h =& \\
	    & \hspace{-3cm} \sum_{i=1}^{N_n} \int_{s_i^n} \ast \bar L_{\boldsymbol \beta}
	    \omega_h \wedge \tilde\eta_h - \sum_{j=1}^{N_{n-1}^{\circ}}
	    \int_{s_{J^{\circ}_j}^{n-1}}
	i_{\boldsymbol\beta}\left(\ast (\omega_h^+-
	\omega_h^-)\wedge \tilde\eta_h^-\right)\,.
        \end{aligned}
    \end{equation}
    Combining these two results for upwind and downwind finite difference we
    get the assertion.
\end{proof}

\begin{remark}
    The careful examination of the proof of Lemma \ref{lemma:standard1} shows
    that the existence of the limit can be shown for quite general
    approximation spaces for differential forms. We used only the continuity
    inside each element.  
\end{remark}

\begin{remark}
    For the sake of completeness we rewrite here the standard convection
    bilinear form for vector proxies $u,v$ or $\ubf, \vbf$ in the case $\Omega
    \subset \mathbb R^3$. $\nbf$ denotes the normal vector of a face.
    Discrete Whitney $0$-forms are continuous. Hence in this case the face
    integrals in \eqref{eq:limit_standard} vanish and we obtain the
    standard bilinear form for scalar convection:
    \begin{equation}
	b_0(\omega_h,\tilde\eta_h)\sim\int_{\Omega}\boldsymbol \beta \cdot \grad
	u(x) \, v(x) dx\,.
    \end{equation}
    For discrete $1$-forms we have both volume and face contributions:
    \begin{multline}
	b_0(\omega_h,\tilde\eta_h)\sim\sum_{i=1}^{N_3} \int_{s_i^3} \grad_h( \boldsymbol \beta \cdot
	\ubf)\cdot \vbf -(\boldsymbol \beta \times \curl \ubf)\cdot \vbf dx \\ -\frac12
	\sum_{j=1}^{N^{\circ}_{2}} \int_{s_{J^\circ_j}^{2}} \boldsymbol \beta \cdot \nbf^+
	(\ubf^+-\ubf^-)\cdot(\vbf^++\vbf^-) dS\,,
    \end{multline}
    and similar for discrete $2$ forms:
    \begin{multline}
	b_0(\omega_h,\tilde\eta_h)\sim\sum_{i=1}^{N_3} \int_{s_i^3} \curl_h(\ubf
	\times \boldsymbol \beta)\cdot \vbf 
	+\boldsymbol \beta \cdot \vbf \, \Div \ubf dx \\ -\frac12
	\sum_{j=1}^{N^{\circ}_{2}} \int_{s_{J^\circ_j}^{2}}\boldsymbol \beta \cdot \nbf^+
	(\ubf^+-\ubf^-)\cdot(\vbf^++\vbf^-) dS\,.
    \end{multline}
    For discrete $3$-forms we get:
    \begin{multline}
	b_0(\omega_h,\tilde \eta_h)\sim\sum_{i=1}^{N_3} \int_{s_i^3} \Div_h(u
	\boldsymbol \beta)v dx  -\frac12
	\sum_{j=1}^{N^{\circ}_{2}} \int_{s_{J^\circ_j}^{2}}\boldsymbol \beta \cdot \nbf^+
	(u^+-u^-)(v^++v^-) dS\,.
    \end{multline}
\end{remark}

\begin{remark}
    For $3$-forms in 3D the limit \eqref{eq:limit_standard} corresponds to
    the non-stabilized discontinuous Galerkin method for the advection
    problem from \cite{BMS04}. In the lowest order case
    this is identical with the scheme for discrete $3$-forms. The limit for the
    upwind finite difference yields the stabilized discontinuous Galerkin
    methods with upwind numerical flux. 
\end{remark}

In practice we will use local quadrature rules to evaluate the volume and
face integrals in $b_0(\omega_h,\tilde \eta_h)$. This will not destroy the
consistency if the quadrature rules are first order consistent.

In the definition of $b_{\Delta\tau}$ from Def.~\ref{def:bd} we may also
appeal to Lemma~\ref{lem:31} and shift the difference quotient onto the
test function. This carries over to the limit and yields two different
Eulerian timestepping methods for problem \eqref{eq:cddf} that,
in analogy to Sect.~\ref{sec:deriv_sLG}, we may call ``direct'' and ``adjoint''.

\begin{definition}
  The direct standard Eulerian timestepping scheme for
  convection-diffusion of differential forms is: given $\omega_h(t_0) \in
  \mathcal W^l(\Omega_h)$ and $\phi \in \mathcal F^l(\Omega)$, find
  $\omega_h(t_k), k>0$ such that:
  \begin{multline}
    \label{def:stEul_cd_dir}
    \int_{\Omega} \ast \frac {\omega_h(t_{k+1})-
      \omega_h(t_k)}{t_{k+1}-t_k}\wedge \tilde\eta_h + \int_{\Omega}
    \ast d \omega_h(t_{k+1})
    \wedge d \tilde\eta \\ + b_0(\omega_h(t_{k+1}),\tilde \eta_h) =
    \int_{\Omega}\ast \phi \wedge \tilde \eta_h\,,\quad \forall \tilde\eta_h \in
    \mathcal W^l(\Omega_h)\,.
  \end{multline}
  The adjoint standard Eulerian timestepping
  is: given $\omega_h(t_0) \in
  \mathcal W^l(\Omega_h)$ and $\phi \in \mathcal F^l(\Omega)$, find
  $\omega_h(t_k), k>0$ such that:
  \begin{multline}
    \label{def:stEul_cd_ad}
    \int_{\Omega} \ast \frac {\omega_h(t_{k+1})-
    \omega_h(t_k)}{t_{k+1}-t_k}\wedge \tilde \eta_h + \int_{\Omega}
    \ast d \omega_h(t_{k+1})
    \wedge d \tilde \eta_h \\ - {b}_0(\tilde \eta_h\,, \omega_h(t_{k+1})) =
    \int_{\Omega} \ast \phi \wedge \tilde \eta,\quad \forall \tilde \eta_h \in \mathcal
    W^l(\Omega_h)\,.
\end{multline}
\end{definition}

\subsection{Upwind scheme}
\label{sec:upwind-scheme}

Now we use the definition \eqref{eq:LieDefinition} of Lie-derivatives as limit of a
difference quotient to define a proper interpolation operator for Lie-derivatives of
discrete forms. Let $s_l^i \in \mathcal S^l_h$ be a subsimplex of the triangulation.
Then both limit values, the limit from the upwind direction
\begin{align}
    \label{eq:upwindLie}
    <L_{\boldsymbol \beta}^{-}\omega_h,s_i^l>:=\lim_{ \Delta \tau \to 0^+}
    \frac{<\omega_h,s_i^l>-<X^*_{t,t-\Delta \tau}\omega_h,s^l_i> }{\Delta
    \tau}\,,
    \quad \omega_h \in \mathcal W^l(\Omega_h)
\end{align}
and the limit from the downwind direction
\begin{align}
    \label{eq:downwindLie}
    <L_{\boldsymbol \beta}^{+}\omega_h,s^l_i>:=\lim_{ \Delta \tau \to 0^+}
    \frac{<X^*_{t,t+\Delta \tau}\omega_h,s^l_i> - <\omega_h,s^l_i>}{\Delta \tau}
    \,,\quad \omega_h \in \mathcal W^1(\Omega_h)
\end{align}
exist for Lipschitz continuous $\boldsymbol \beta$. The continuity of
traces of discrete differential forms ensures that integrals
are well-defined. The existence of upwind and downwind limits then follows
from the existence of the limits for locally linear $\boldsymbol \beta$, which
in turn can be calculated explicitly by distinguishing different geometric
situations. We skip the tedious details.

Observe that in general
\begin{align}
    <L_{\boldsymbol \beta}^{+}\omega_h,s^l_i>\neq <L_{\boldsymbol
    \beta}^{-}\omega_h,s^l_i>\,,
\end{align}
e.g. for $0$-forms a short calculation shows
\begin{align}
    <L_{\boldsymbol \beta}^{-}\omega_h,s_i^0>=\boldsymbol \beta(s_i^0) \cdot \grad
    \omega_{h}^-(s_i^0)\;,
\end{align}
and 
\begin{align}
    <L_{\boldsymbol \beta}^{+}\omega_h,s_i^0>=\boldsymbol \beta(s_i^0) \cdot \grad
    \omega_{h}^+(s_i^0)\;,
\end{align}
with upwind and downwind traces $\omega_{h}^-$ and $\omega_{h}^+$.
Only for cases where $\omega_h$ is continuous in a neighbourhood of subsimplex
$s_i^l$ upwind and downwind limit agree.

Since \eqref{eq:upwindLie} and \eqref{eq:downwindLie} correspond to integral
evaluation of the usual interpolation operators, we could either use the upwind or the
downwind limit to define discrete Lie-derivatives. We will use the upwind limit for
the direct formulation and the downwind limit for the adjoint formulation.  This is
motivated by taking the cue from semi-Lagrangian schemes.

\begin{definition}
    Let $\boldsymbol \beta$ be Lipschitz continuous and $\omega_h \in
    \mathcal W^l(\Omega_h)$. Further let $(s^l_i)_{i=1}^{N_l}$ be
    the subsimplices corresponding to the degrees of freedom of basis
    functions $(b_{i}^l)_{i=1}^{N_l}$ then 
    \begin{align}
	L^{-}_{\boldsymbol \beta,h} \omega_h:= \sum_i
	<L_{\boldsymbol \beta}^{-} \omega,s_i^l> b_i^l
	\label{eq:interpol_ufd}
    \end{align}
    is the upwind interpolated discrete Lie-derivative.
    Analogously the downwind interpolated Lie-derivative is
    \begin{align}
	L^{+}_{\boldsymbol \beta,h} \omega_h:= \sum_i
	<L_{\boldsymbol \beta}^{+} \omega,s_i^l> b_i^l\,.
    \end{align}
\end{definition}
It follows from the Cartan formula \eqref{eq:Cartan} that 
\begin{align}
    <L_{\boldsymbol \beta}^- \omega_h,s_i^l>=<d (i_{\boldsymbol \beta}^-
    \omega_h) + i_{\boldsymbol \beta}^- (d\omega_h),s_i^l >\,,
    \label{}
\end{align}
with
\begin{align}
    <i_{\boldsymbol \beta}^- \omega_h ,s_i^l>:=\lim_{\Delta \tau \to 0^+}
    <\omega_h,Ext_{t,t-\Delta \tau}(\boldsymbol \beta, s_i^l)>\,.
    \label{}
\end{align}

\begin{remark}
The exterior derivatives can be evaluated exactly for discrete differential forms.
Thus, the contraction operators are the only approximations. But for locally constant
discrete forms these approximations are exact and one can prove consistency
\begin{align}
    |\int_\Omega \ast L_{\boldsymbol \beta} \wedge \tilde \eta -\int_{\Omega}
    \ast L^{-}_{\boldsymbol \beta, h} \Pi^l \omega \wedge \Pi^l \tilde \eta |
    = O(h)\,,
    \quad \forall \omega, \tilde \eta \in \mathcal F^l(\Omega)\,,
\end{align}
once one has interpolation estimates for discrete differential forms. This is
outside the scope of this paper. 
\end{remark}

\begin{figure}[h]
    \begin{center}
	\psfrag{1}{$s_1^0$}
	\psfrag{2}{$s_2^0$}
	\psfrag{3}{$s_3^0$}
	\psfrag{4}{$s_4^0$}
	\psfrag{5}{$s_5^0$}
	\psfrag{6}{$s_6^0$}
	\psfrag{7}{$s_1^1$}
	\psfrag{8}{$s_2^1$}
	\includegraphics[width=0.5\linewidth]{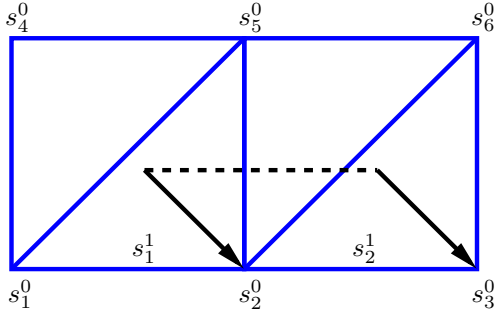}
    \end{center}
    \caption{Example that we can not change the order of integration and
    limit.}
    \label{fig:limit_example}
\end{figure}

\begin{remark}
  We finally want to stress that it is important to treat the limits
  \eqref{eq:upwindLie} and \eqref{eq:downwindLie} for $l>0$ as limits of co-chains
  rather than integrals with certain integrands that are defined as pointwise limits.
  First, these pointwise limits are not well-defined since discrete forms are in
  general discontinuous across element boundaries. Second, even for special cases
  where the pointwise limit $\lim_{\Delta \tau \to 0} \frac{\omega_h-X^*_{t,t-\Delta
      \tau}\omega_h}{\Delta \tau}(x)$ exist we can have
    \begin{multline}
	\lim_{ \Delta \tau \to 0^+}
	<\frac{\omega_h-X^*_{t,t-\Delta \tau}\omega_h }{\Delta
	\tau} ,s^l_i>\neq \\ <\lim_{\Delta \tau \to 0^+} \frac{\omega_h-X^*_{t,t-\Delta
	\tau}\omega_h}{\Delta \tau},s^l_i>\,,
	\quad \omega_h \in \mathcal W^l(\Omega_h)\,.
	\label{}
    \end{multline}
    To understand this, we may consider the 2D example depicted in
    Fig.~\ref{fig:limit_example} with constant $\boldsymbol \beta$. We take $\omega_h
    \in \mathcal W^1$ such that $<\omega_h,s_{i}^1>=0$ for all edges $s_i^1\in
    \mathcal S^1$ except the edge $s_1^1$ between vertices $s_1^0$ and $s_2^0$ and
    calculate the projection of the Lie-derivative onto edge $s_2^1$ between vertices
    $s_2^0$ and $s_3^0$.  Then it is clear that
    \begin{align}
	\lim_{\Delta \tau \to 0^+} \frac{\omega_h-X^*_{t,t-\Delta
	\tau}\omega_h}{\Delta \tau}(x)& =0\,, \quad x\in s_2^1\cup\{s_3^0\}
    \end{align}
    and therefor
    \begin{align}
	<\lim_{\Delta \tau \to 0^+} \frac{\omega_h-X^*_{t,t-\Delta
	\tau}\omega_h}{\Delta \tau},s_2^1> = 0\,.
    \end{align}
    On the other hand, one can easily show that:
    \begin{equation}
	\begin{aligned}
	    <\omega_h-X^*_{t,t-\Delta \tau}\omega_h 
	    ,s^1_2>\\
	    & \hspace{-3cm} =<\omega_h ,X_{t,t-\Delta \tau} s^1_2> \\
	    & \hspace{-3cm} = \lambda_1(s_2^0-\Delta \tau \boldsymbol \beta)
	    \lambda_2(s_2^0-\frac{\Delta \tau}{\sqrt{2}}
	    \boldsymbol \beta)- \lambda_1(s_2^0-\frac{\Delta \tau}{\sqrt{2}} \boldsymbol \beta)
	    \lambda_2(s_2^0-\Delta \tau \boldsymbol \beta) \\
	    &\hspace{-3cm} = \Delta \tau (\frac{1}{\sqrt
	    2}\grad \lambda_1 \cdot \boldsymbol \beta-\grad \lambda_1 \cdot
	    \boldsymbol \beta)+ O(\Delta \tau ^2)\,, 
	    \label{}
	\end{aligned}
    \end{equation}
    hence
    \begin{align}
	\lim_{ \Delta \tau \to 0^+}
	<\frac{\omega_h-X^*_{t,t-\Delta \tau}\omega_h }{\Delta
	\tau} ,s^1_2>=\frac{1-\sqrt{2}}{\sqrt{2}} \grad \lambda_1 \cdot
	\boldsymbol \beta\neq 0\,.
	\label{}
    \end{align}
    Since an explicit calculation of the upwind limits
    \eqref{eq:upwindLie} might be expensive for arbitrary
    $\boldsymbol \beta$, we have to find consistent
    approximations. The preceding discussion shows that this must be done very
    carefully. Replacing the integration simply with some quadrature will not
    yield a consistent approximation. 
\end{remark}

We combine the upwind discretization with an implicit Euler method and state
the implicit upwind methods for convection-diffusion for
discrete differential forms. 

\begin{definition}
    The implicit direct upwind Eulerian timestepping scheme for
    the variational formulation of convection-diffusion of differential
    forms is: Given $\omega_h(t_0) \in
    \mathcal W^l(\Omega_h)$ and $\phi \in \mathcal F^l(\Omega)$, find
    $\omega_h(t_k), k>0$ such that:
    \begin{multline}
	\label{def:iuEul_cd_dir}
	\int_{\Omega} \ast \frac {\omega_h(t_{k+1})-
	\omega_h(t_k)}{t_{k+1}-t_k}\wedge \tilde \eta + \varepsilon \int_{\Omega}
	\ast d \omega_h(t_{k+1})
	\wedge d \tilde \eta \\ + \int_{\Omega} \ast L^-_{\boldsymbol \beta,h
	}\omega_h(t_{k+1}) \wedge \tilde \eta =
	\int_{\Omega} \ast \phi \wedge \tilde \eta\,,\quad \forall \tilde \eta \in
	\mathcal W^l(\Omega_h)\,.
    \end{multline}
    The adjoint upwind Eulerian timestepping scheme for
    the adjoint variational formulation of convection-diffusion of differential
    forms is: Given $\omega_h(t_0) \in
    \mathcal W^l(\Omega_h)$ and $\phi \in \mathcal F^l(\Omega)$, find
    $\omega_h(t_k), k>0$ such that:
    \begin{multline}
	\label{def:iuEul_cd_ad}
	\int_{\Omega} \ast \frac {\omega_h(t_{k+1})-
	\omega_h(t_k)}{t_{k+1}-t_k}\wedge \tilde \eta +\varepsilon \int_{\Omega}
	\ast d \omega_h(t_{k+1})
	\wedge d \tilde \eta \\- \int_{\Omega} \ast \omega_h(t_{k+1}) \wedge
	L^{+}_{\boldsymbol \beta,h}\tilde \eta =
	\int_{\Omega} \ast \phi \wedge \tilde \eta\,, \quad \forall \tilde \eta \in
	\mathcal W^l(\Omega_h)\,.
    \end{multline}
\end{definition}

In contrast to the semi-Lagrangian schemes \eqref{def:SL_cd_dir} and
\eqref{def:SL_cd_ad}, in \eqref{def:iuEul_cd_dir} and \eqref{def:iuEul_cd_ad} we face
a non-symmetric algebraic system in each time step, whose iterative solution can be
challenging. We therefor consider a second kind of upwind Eulerian scheme, that
treats the convection in an explicit fashion.

\begin{definition}
    The direct semi-implicit upwind Eulerian timestepping scheme for
    the variational formulation of convection-diffusion of differential
    forms is: Given $\omega_h(t_0) \in
    \mathcal W^l(\Omega_h)$ and $\phi \in \mathcal F^l(\Omega)$, find
    $\omega_h(t_k), k>0$ such that:
    \begin{multline}
	\label{def:siuEul_cd_dir}
	\int_{\Omega} \ast \frac {\omega_h(t_{k+1})-
	\omega_h(t_k)}{t_{k+1}-t_k}\wedge \tilde \eta  + \varepsilon \int_{\Omega}
	\ast d \omega_h(t_{k+1})
	\wedge d \tilde \eta \\ + \int_{\Omega} \ast L^-_{\boldsymbol \beta,h
	}\omega_h(t_{k}) \wedge \tilde \eta =
	\int_{\Omega} \ast \phi \wedge \tilde \eta\,, \quad \forall \tilde \eta \in
	\mathcal W^l(\Omega_h)\,.
    \end{multline}
    The adjoint upwind Eulerian timestepping scheme for
    the variational formulation of convection-diffusion of differential
    forms is: Given $\omega_h(t_0) \in
    \mathcal W^l(\Omega_h)$ and $\phi \in \mathcal F^l(\Omega)$, find
    $\omega_h(t_k), k>0$ such that:
    \begin{multline}
	\label{def:siuEul_cd_ad}
	\int_{\Omega} \ast \frac {\omega_h(t_{k+1})-
	\omega_h(t_k)}{t_{k+1}-t_k}\wedge \tilde \eta  +\varepsilon \int_{\Omega}
	\ast d \omega_h(t_{k+1})
	\wedge d \tilde \eta \\ - \int_{\Omega} \ast \omega_h(t_{k}) \wedge 
	L^{+}_{\boldsymbol \beta,h} \tilde \eta =
	\int_{\Omega} \ast \phi \wedge \tilde \eta\,, \quad \forall \tilde \eta \in
	\mathcal W^l(\Omega_h)\,.
    \end{multline}
\end{definition}

These semi-implicit schemes do not coincide with the usual method of lines approach.
Nevertheless there are two heuristic arguments that could justify the semi-implicit
schemes.  Note first that for discrete $0$-forms the Semi-Lagrangian and the
semi-implicit Eulerian scheme are identical, if we use a time step length $\Delta
t=O(h)$ and explicit Euler for solving the characteristic equation
\eqref{eqn:charac_ode}.  
The direct Semi-Lagrangian scheme \eqref{def:SL_cd_dir} for discrete $0$-forms was,
written here for vector proxies $u$ and $v$ for $\omega_h$ and $\eta_h \in \mathcal
W^0(\Omega_h)$: Find $u$ such that
\begin{multline*}
    \int_{\Omega} \frac{u(t_{k+1},x)-\Pi^0
    u(t_{k},\bar{X}_{t_{k+1},t_{k}}(x))}{t_{k+1}-t_{k}} v(x) dx\\+\varepsilon
    \int_{\Omega}\grad u(t_{k+1},x) \cdot \grad v(x) dx = \int_{\omega} f(x)
    v(x) dx\,, \quad
    \forall v \in \mathcal W^0(\Omega_h)
\end{multline*}
while the semi-implicit Eulerian scheme was: Find $u$ such that 
\begin{multline*}
    \int_{\Omega} \frac{u(t_{k+1},x)-
    u(t_{k},x)}{t_{k+1}-t_{k}} v(x) dx+\varepsilon
    \int_{\Omega}\grad u(t_{k+1},x) \cdot \grad v(x) dx \\+ \int_{\Omega} L^-_{\boldsymbol
    \beta,h}u(t_{k},x) v(x)
    dx  = \int_{\Omega} f(x) v(x) dx\,, \quad
    \forall v \in \mathcal W^0(\Omega)\,,
\end{multline*}
with 
\begin{eqnarray*}
    L_{\boldsymbol \beta,h}^-u(t_{k},x)=\sum_{i=1}^{N_0} \boldsymbol
    \beta(s_{i}^0) \cdot
    \grad u(t_{k},s_i^0)_{|T_i} \lambda_i(x)\,.
\end{eqnarray*}
On the other hand $\bar{X}_{t_{k+1},t_{k}}(s_i^0)= s_i^0- (t_{k+1}-t_{k})\boldsymbol \beta
(s_i^0)$ and by the time step condition $\Delta t = O(h)$:
\begin{multline*}
    \Pi^0 u(t_{k},\bar{X}_{t_{k+1},t_{k}}(x))=\sum_{i=1}^{N_0}
    u(t_{k},s_i^0)\lambda_i(x)\\-(t_{k+1}-t_{k}) \boldsymbol \beta (s_i^0) \cdot \grad
    u(t_{k},s_i^0)_{|T_i} \lambda_i(x)\,,
\end{multline*}
hence we have the identity:
\begin{eqnarray*}
    \Pi^0
    u(t_{k},\bar{X}_{t_{k+1},t_{k}}(x))=u(t_{k},x)-(t_{k+1}-t_{k})L^-_{\boldsymbol
    \beta,h} u(t_{k})\,.
\end{eqnarray*}
The second argument, that supports the definition of semi-implicit schemes
applies to forms of any degree and follows from the identity:
\begin{eqnarray*}
    X_{t,t'}\omega&=&\omega-(X_{t,t'}\omega-\omega)\\
    &=& \omega-(t'-t)L_{\boldsymbol \beta}\omega+O( (t'-t)^2),\quad \omega \in
    \mathcal F^l(\Omega)\,.
\end{eqnarray*}
Hence we can derive the semi-implicit schemes \eqref{def:iuEul_cd_dir} and
\eqref{def:iuEul_cd_ad} from the Semi-Lagrangian schemes \eqref{def:SL_cd_dir} and
\eqref{def:SL_cd_ad} in replacing e.g. 
\begin{eqnarray*}
    \Pi^l\bar{X}_{t_{k+1},t_k}\omega_h(t_k)\approx
    \omega_h(t_k)-(t_k-t_{k+1})L^-_{\boldsymbol \beta,h}\omega_h(t_k)\,.   
\end{eqnarray*}

\begin{remark}
  The analysis of a semi-implicit discontinuous Galerkin scheme for scalar non-linear
  convection-diffusion can be found in \cite{DFH07}. For the linear case this scheme is
  modulo certain quadrature rules identical with the semi-implicit Euler scheme for
  $3$-forms.
\end{remark}

\begin{remark}
  There is also a very close link between both the semi-Lagrangian and
  semi-implicit Eulerian methods and arbitrary Lagrangian-Eulerian methods
  (ALE) \cite[p. 322]{Benson}. 
  ALE methods are operator splitting methods, that treat the
  diffusion part in a Lagrangian and the convection part in an Eulerian
  fashion. After a certain number of Lagrangian iterations steps an Eulerian
  iteration step maps the mesh function on the distorted mesh back to the
  initial mesh. The semi-Lagrangian and semi-implicit
  Eulerian methods apply this mapping in each time step.
  In many problems these remap operations should preserve certain properties of the mesh
  function, e.g. volume or closedness. It is the discrete differential form
  approach, that naturally guaranties such properties for the methods presented
  here and for the ALE method in \cite{RWW07}.   
\end{remark}

\section{Application: Magnetic convection-diffusion}
\label{sec:MHDmodel}

In this section we will derive two different convection-diffusion equations for the
electromagnetic part of magnetohydrodynamics (MHD) models. The first one is an
equation for the magnetic field $h$ and requires the adjoint methods while the second
one for the magnetic vector potential $a$ can be solved via the direct methods.

Commonly in MHD applications, one neglects the displacement current. This reduced
model, called eddy current model, is a system of equations for the magnetic field $h
\in {\mathcal F}^1(\Omega)$, the electric field $e \in {\mathcal F}^1(\Omega)$, the
magnetic field density $b \in {\mathcal F}^2(\Omega)$, the current density $j \in
{\mathcal F}^2(\Omega)$ and external source $f \in {\mathcal F}^2(\Omega)$:
\begin{align}
    d e& =-\partial_t b  &\text{in $\Omega$} \label{eq:fala}\,, \\
    d h& = j + f & \text{in $\Omega$}\,, \\ 
    j&=\ast_{\sigma}(e-i_{\boldsymbol \beta} b)&\text{in $\Omega$}\,, \\
    \ast_{\mu}h&=b & \text{in $\Omega$}\,. 
\end{align}
A uniformly positive conductivity $\sigma$ will be assumed in the sequel.

To rewrite the eddy current model in terms of a material derivative we substitute
$e=\tilde{e}+i_{\boldsymbol \beta}b$ and add $-i_{\boldsymbol \beta}db$ to 
Faraday's law (\ref{eq:fala}). This leaves the solution unchanged, since $db=0$.
Hence we end up with the system:
\begin{align}
    d \tilde{e}& =-D_{\boldsymbol \beta} b &\text{in $\Omega$}\,, \label{eq:fl} \\
    d h& = j+f& \text{in $\Omega$}\,, \label{eq:am} \\ 
    j&=\ast_{\sigma}\tilde{e}& \text{in $\Omega$}\,, \label{eq:Mat_sigma}\\
    \ast_{\mu}h&=b& \text{in $\Omega$}\,.\label{eq:Mat_mu}
\end{align}
Next we state the two different variational formulations, \textit{cf.} \cite[Sect.~2.3]{HIP02}.

\subsection{$h$-based variational formulation}
\label{section:hform}
For simplicity we impose homogeneous electric boundary conditions on $\partial
\Omega$. Testing (\ref{eq:fl}) with $h'\in \mathcal F^1(\Omega)$, integration by parts
yields:
\begin{align}
    \int_\Omega \tilde{e} \wedge d h'+ \int_\Omega D_{\boldsymbol \beta} b
    \wedge h'=0\,. 
\end{align}
We eliminate $\tilde{e}$ using Ohm's law and end up with the following
variational formulation:
Seek $h \in \mathcal F^1(\Omega), b \in \mathcal F^2(\Omega)$ such that:
\begin{align}
    \int_\Omega D_{\boldsymbol \beta} b \wedge h'+\int_\Omega \ast_{\sigma^{-1}}(dh-f)\wedge
    dh'& =0 \,,\quad \forall h' \in
    \mathcal F^1(\Omega)\,, \\
    \int_\Omega \ast_{\mu}h\wedge h''-\int_\Omega b \wedge h'' &= 0 \,,\quad
    \forall h'' \in \mathcal
    F^1(\Omega)\,.
\end{align}
In order to eliminate $b$ as well we use Corollary \ref{cor:dualLie} and get 
\begin{multline}
    \int_\Omega (\ast_{\mu}\partial_t h \wedge h' - \ast_{\mu} h \wedge   
    L_{\boldsymbol \beta} h') + \int_\Omega \ast_{\sigma^{-1}}(dh-f) \wedge dh' 
    = 0\,, \quad \forall h' \in
    \mathcal F^1(\Omega)
\end{multline}
or, equivalently
\begin{multline}
    \int_\Omega \lim_{\Delta t} \frac{\ast_{\mu}h(t) \wedge h'-
    \ast_{\mu} h(t-\Delta t) \wedge
    X^{*}_{t-\Delta t,t} h' }{\Delta t}+   
    \int_\Omega \ast_{\sigma^{-1}}(dh-f) \wedge dh' 
    = 0\,, \quad \forall h' \in
    \mathcal F^1(\Omega)\,.
\end{multline}
This variational formulation can be approximated using a Semi-Lagrangian
framework \eqref{def:SL_cd_ad} yielding an algebraic system:
\begin{align}
    \label{eq:h_timestepping}
    (\mathbf M_\mu +\Delta t\mathbf C_{\sigma^{-1}})\mathbf h^t =
    ({\mathbf P}_{t-\Delta,t}^1)^T\mathbf
    M_\mu \mathbf h^{t-\Delta t} +\Delta t \mathbf f\,,
\end{align}
where $\Mbf_\mu$ and $\Cbf_{\sigma^{-1}}$ are the mass matrices
$\int_{\Omega}\ast_{\mu}\omega_h \wedge \tilde \eta_h$ and
$\int_{\Omega}\ast_{\sigma^{-1}} d \omega_h \wedge d \tilde \eta_h$, $\omega_h,
\tilde \eta_h \in \mathcal W^1(\Omega_h)$.  The evaluation of the right hand side
requires the calculation of the matrix ${\mathbf P}^1_{\tau,t}$ from
\eqref{eq:Pmat1F}.  Lemma (\ref{lemma:commute}) shows that the solution of
(\ref{eq:h_timestepping}) is weakly closed if the initial data is weakly closed and
$f=0$.  Alternatively, we may use one of the implicit Eulerian schemes
\eqref{def:stEul_cd_ad} or \eqref{def:iuEul_cd_ad}
\begin{align}
    \label{eq:h_timestepping_Eulerian1}
    (\mathbf M_\mu +\Delta t\mathbf C_{\sigma^{-1}}-\Delta t \mathbf
    L_a )\mathbf h^t = \mathbf M_\mu \mathbf h^{t-\Delta t} + \Delta t \mathbf
    f\,,
\end{align}
where $\mathbf L_a$ is the stiffness matrix for either the standard
\eqref{def:stEul_cd_ad} or the
downwind \eqref{def:iuEul_cd_ad} discretization.
While these schemes will not preserve the weakly closed property, the
semi-implicit upwind Eulerian scheme \eqref{def:siuEul_cd_ad} with algebraic system
\begin{align}
    \label{eq:h_timestepping_Eulerian2}
    (\mathbf M_\mu +\Delta t \mathbf C_{\sigma^{-1}}
    )\mathbf h^t = (\mathbf M_\mu +\Delta t \mathbf
    L_a) \mathbf h^{t-\Delta t} + \Delta t \mathbf f
\end{align}
does, due to the interpolation based construction. 

\subsection{$a$-based variational formulation}
Here, for simplicity, we assume homogeneous magnetic boundary conditions on
$\partial \Omega$. The ansatz $\tilde{e}=-D_{\boldsymbol \beta}a$ and $b=da$ solves Faraday's
law (\ref{eq:fl}). We then multiply Ampere's law with a test form $a' \in \mathcal
D^1$ and integrate by parts:
\begin{align}
    \int_\Omega h \wedge d a' - \int_\Omega (j+f) \wedge a' =0\,, \quad \forall a'
    \in \mathcal F^1(\Omega)\,.
\end{align}
The material laws (\ref{eq:Mat_sigma}) and ($\ref{eq:Mat_mu}$) eliminate
$j$ and $h$ and we get the $a$-based variational formulation:
Find $a \in \mathcal F^1(\Omega)$ such that:
\begin{align}
    \int_\Omega \ast_{\mu^{-1}}da \wedge da'+\int_\Omega
    \ast_{\sigma}D_{\boldsymbol \beta} a
    \wedge a'= \int_\Omega f \wedge a'\,, \quad  \forall a' \in
    \mathcal F^1(\Omega)\,.
\end{align}
The algebraic system for the direct Semi-Lagrangian scheme
\eqref{def:SL_cd_dir} then reads as 
\begin{align}
    ( \mathbf M_\sigma+\Delta t \mathbf C_{\mu^{-1}} )\mathbf a^t =  \mathbf
    M_\sigma {\mathbf P}_{t,t-\Delta t}^1\mathbf a^{t-\Delta t}+ \Delta t
    \mathbf f\,,
\end{align}
while the systems for the implicit and semi-implicit Eulerian schemes are:
\begin{align}
    \label{eq:a_timestepping_Eulerian1}
    (\mathbf M_\sigma + \Delta t \mathbf C_{\mu^{-1}}+\Delta t 
    \mathbf L_d)\mathbf a^t = \mathbf M_\sigma
    \mathbf a^{t-\Delta t} + \Delta t \mathbf f
\end{align}
and
\begin{align}
    \label{eq:a_timestepping_Eulerian2}
    (\mathbf M_\sigma +\Delta t \mathbf C_{\mu^{-1}})\mathbf a^t = (\mathbf
    M_\sigma - \Delta t \mathbf L_d) \mathbf a^{t-\Delta t}+\Delta \mathbf f\,.
\end{align}
Here $\Mbf_{\sigma}$ and $\Cbf_{\mu}^{-1}$ are the stiffness matrices for
$\int_{\Omega}\ast_{\sigma}\omega_h \wedge \tilde \eta_h$  and
$\int_{\Omega}\ast_{\mu^{-1}} d \omega_h \wedge d \tilde \eta_h$,
$\omega_h, \tilde \eta_h \in \mathcal W^1(\Omega_h)$.
$\mathbf L_d$ corresponds either to the standard \eqref{def:stEul_cd_dir} or the upwind
Eulerian \eqref{def:iuEul_cd_dir} discretization of the convection part. 

\section{Numerical experiments}
We finally we present a few numerical examples that illustrate the performance
of the methods derived here. We take $\Omega \subset \mathbb R^2$ and look at
\begin{gather}
  \begin{array}[c]{rcll}
    \partial_t \ast \omega(t) - \varepsilon (-1)^{l} 
    d\ast d \omega(t) +  
    L_{\boldsymbol \beta} \ast\omega (t) &=& \varphi \quad &
    \text{in }\Omega\subset\mathbb{R}^{2}\;,\\
    \imath^{\ast}\omega &=& 0 & \text{on }\partial\Omega\;,\\
    \omega(0) &=& \omega_{0} & \text{in } \Omega \;
  \end{array}
\end{gather}
for 1-forms $\omega \in \mathcal F^1(\Omega)$. In vector proxy notion with
$\ubf:=\operatorname{v.p.}(\omega)$ this reads
\begin{gather}
    \label{eq:problem1}
    \begin{array}[c]{rcll}
	\partial_{t}\ubf + \varepsilon \curl \text{curl} \ubf + {\boldsymbol
	\beta}(\Div \ubf) + \curl (\ubf \times  {\boldsymbol \beta}) &=& \fbf & \text{in
	}\Omega\;,\\
	\ubf\cdot\nbf &=& 0 & \text{on }\partial\Omega\;,\\
	\ubf(0) &=& \ubf_{0} & \text{in } \Omega\;,
    \end{array}
\end{gather}
with $\curl := \Rbf\grad$, $\text{curl}:=\Div \Rbf$ and $\ubf \times \boldsymbol
\beta := \ubf^T \Rbf\boldsymbol \beta$, $\Rbf = \begin{pmatrix} 0 & -1 \\ 1 & 0
\end{pmatrix}$.  We approximate $\omega$ by discrete $1$-forms $\omega_h$ on a
triangular mesh. We study the adjoint Semi-Lagrangian and Eulerian methods
\eqref{def:SL_cd_ad}, \eqref{def:stEul_cd_ad}, \eqref{def:iuEul_cd_ad} and
\eqref{def:siuEul_cd_ad} for this boundary value problem.

\begin{figure}
    \begin{center}
	\includegraphics[width=0.75\linewidth]{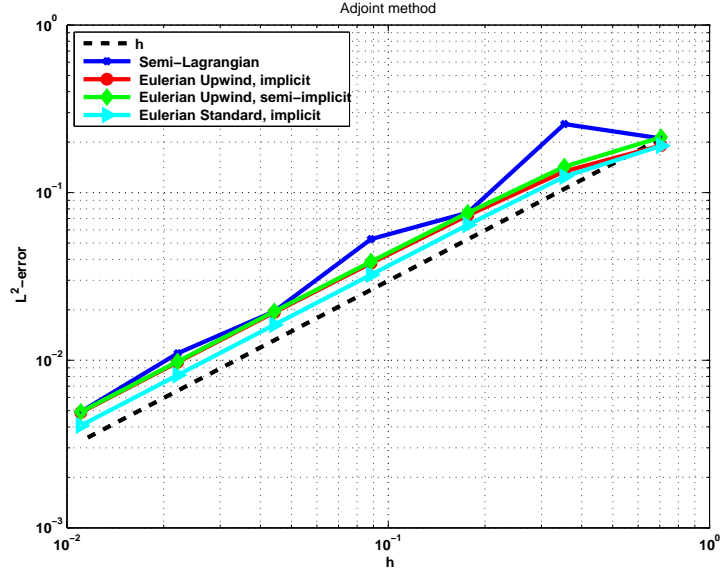}
    \end{center}
    \caption{Experiment I: $L^2$-error at $t=0.25$ for the problem with divergence free
    velocity and small
    CFL-number $CFL=0.1$, $\varepsilon=1$.}
    \label{fig:convplot1}
\end{figure}

\begin{figure}
    \begin{center}
	\includegraphics[width=0.75\linewidth]{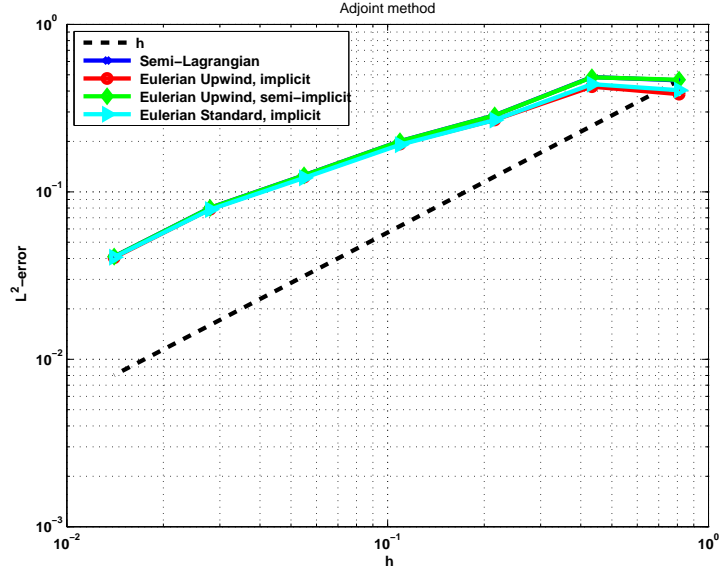}
    \end{center}
    \caption{Experiment I: $L^2$-error at $t=0.25$ for the problem with divergence free velocity and larger
    CFL-number $CFL=0.8$, $\varepsilon=1$.}
    \label{fig:convplot2}
\end{figure}

\textbf{Experiment I.} In the first experiment we take in problem \eqref{eq:problem1} the domain
$\Omega=[-1,1]^2$, the divergence free
velocity
\begin{eqnarray}
    \label{eq:velo1}
    \boldsymbol \beta =
    \begin{pmatrix}
	(1-x^2)^2(y-y^3) \\
	-(1-y^2)^2(x-x^3)
    \end{pmatrix}
\end{eqnarray}
and 
\begin{eqnarray*}
  \ubf =\cos(2 \pi t) 
    \begin{pmatrix} 
	\sin(\pi x)\sin(\pi y) \\  (1-x^2) (1-y^2) 
    \end{pmatrix}.
\end{eqnarray*}
Then with $\varepsilon=1$ neither the convection nor the diffusion part
dominates and we expect for all three schemes first order convergence in the mesh size $h$, if
the time step $\Delta t: = t-\tau$ is of order $h$. The convergence plots in
figures \ref{fig:convplot1} and \ref{fig:convplot2} show the $L^2$-error
at time $t=0.25$ and confirm this for
different ratios $CFL:=\frac{\Delta t \|\boldsymbol \beta\|_{\infty}}{h}$.

\begin{figure}
    \begin{center}
	\includegraphics[width=0.75\linewidth]{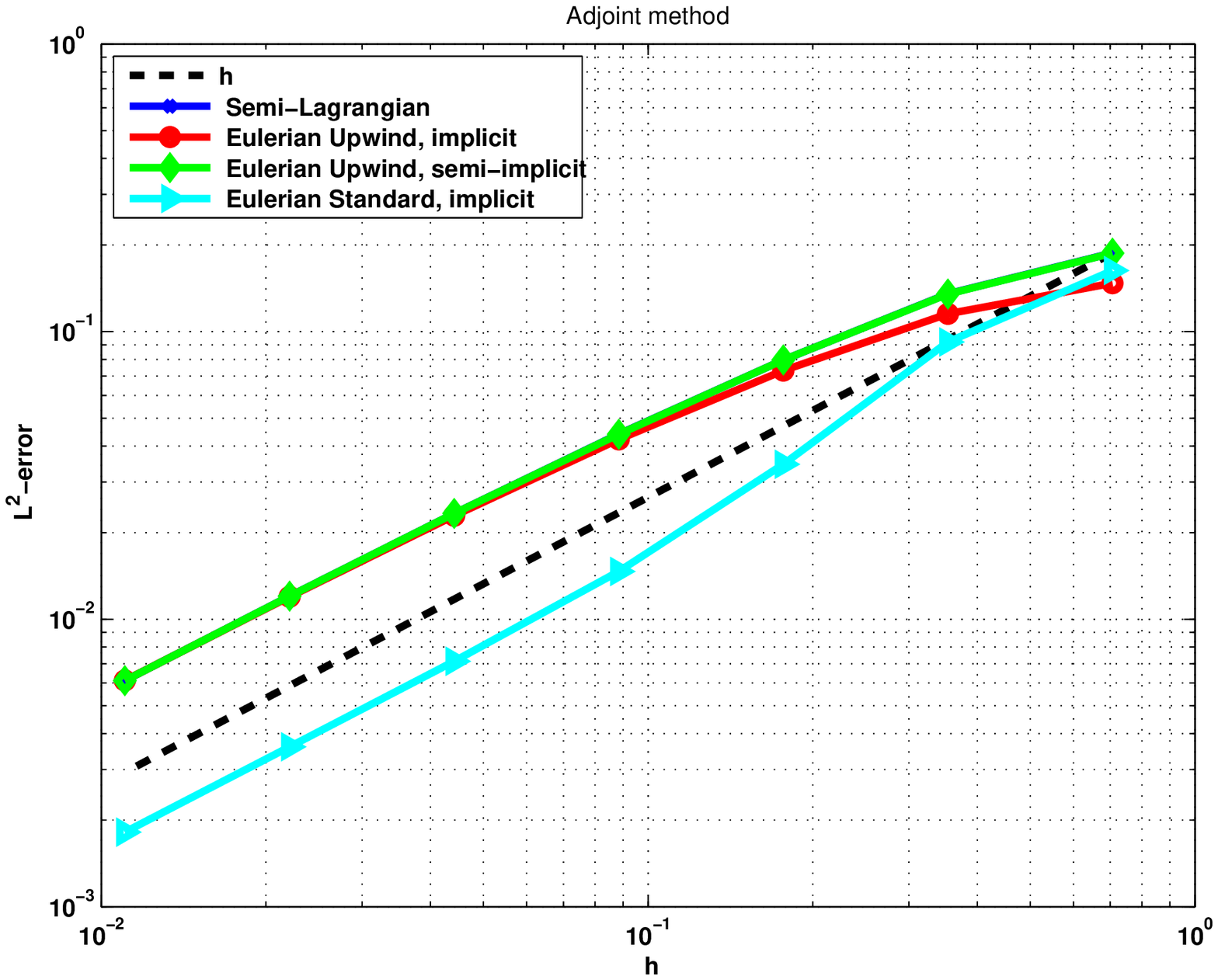}
    \end{center}
    \caption{Experiment II: $L^2$-error at $t=0.25$ for the problem with non-divergence free velocity and small
    CFL-number $CFL=0.1$, $\varepsilon=1$.}
    \label{fig:convplot3}
\end{figure}

\begin{figure}
    \begin{center}
	\includegraphics[width=0.75\linewidth]{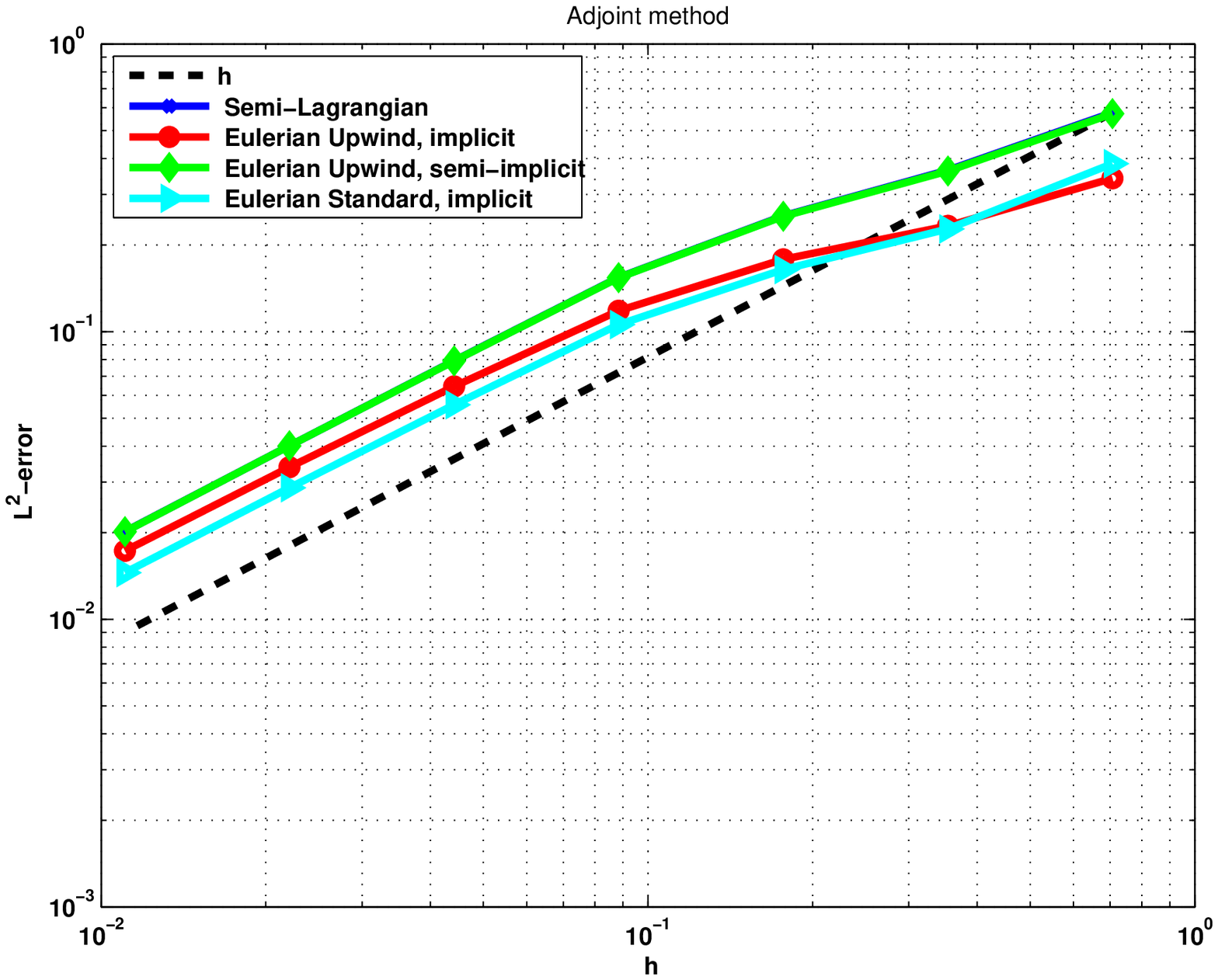}
    \end{center}
    \caption{Experiment II: $L^2$-error at $t=0.25$ for the problem with non-divergence free velocity and
    larger CFL-number $CFL=0.8$,$\varepsilon=1$.}
    \label{fig:convplot4}
\end{figure}

\textbf{Experiment II.} We consider the non-divergence free velocity field
\begin{eqnarray*}
    \boldsymbol \beta = 
    \begin{pmatrix}
	sin(\pi x) (1-y^2) \\
	sin(\pi y) (1-x^2)
    \end{pmatrix}\;.
\end{eqnarray*}
Again, we obtain first order convergence (see figs \ref{fig:convplot3}) and
\ref{fig:convplot4}).

\begin{figure}
    \begin{center}
	\includegraphics[width=0.75\linewidth]{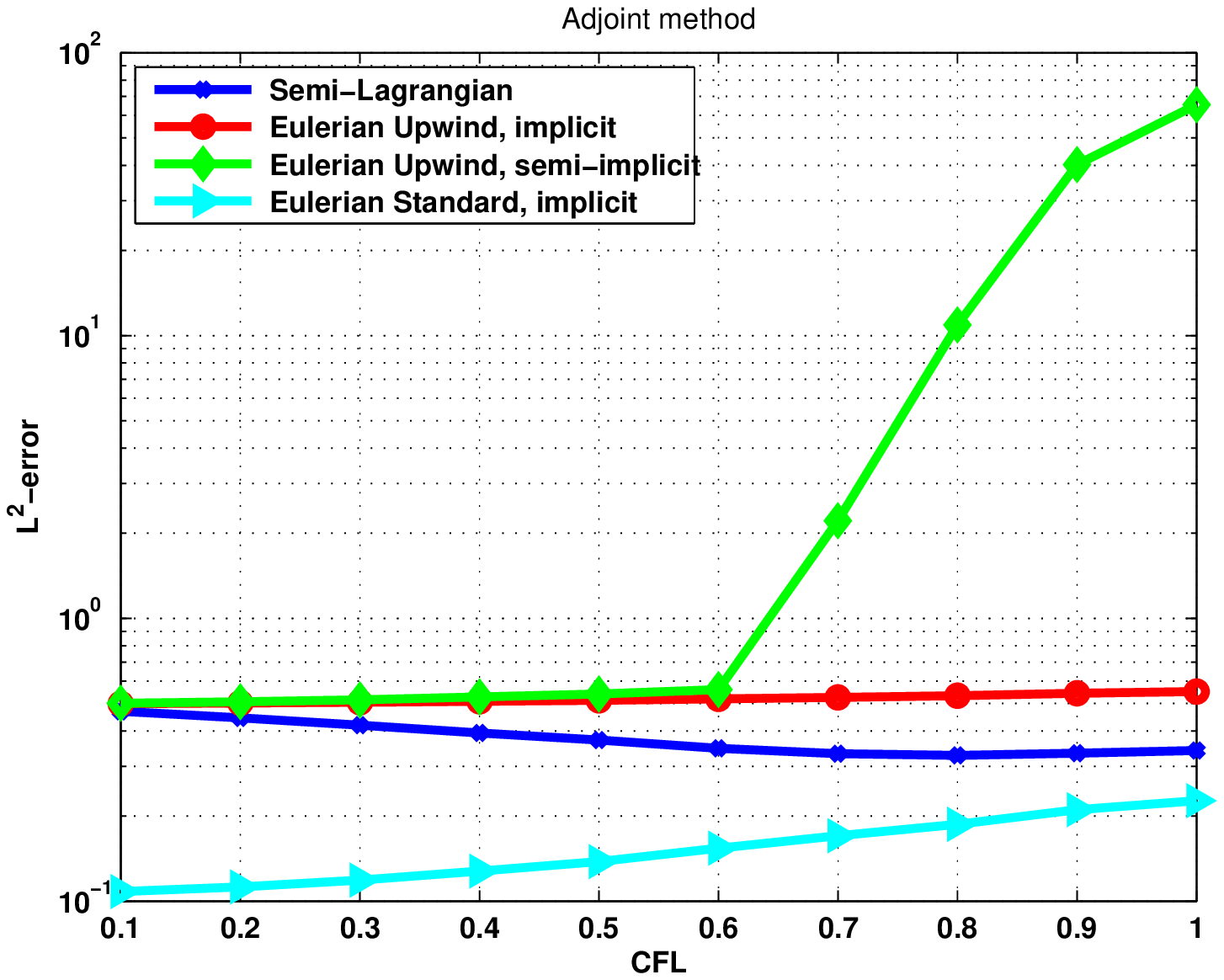}
    \end{center}
    \caption{Experiment III: $L^2$-error for different CFL-numbers at $t=0.5$ for evolution in
    fixed triangulation with mesh size $h_m=0.11$ with
    $\varepsilon=10^{-3}$.}
    \label{fig:stab1}
\end{figure}

\begin{figure}
    \begin{center}
	\includegraphics[width=0.75\linewidth]{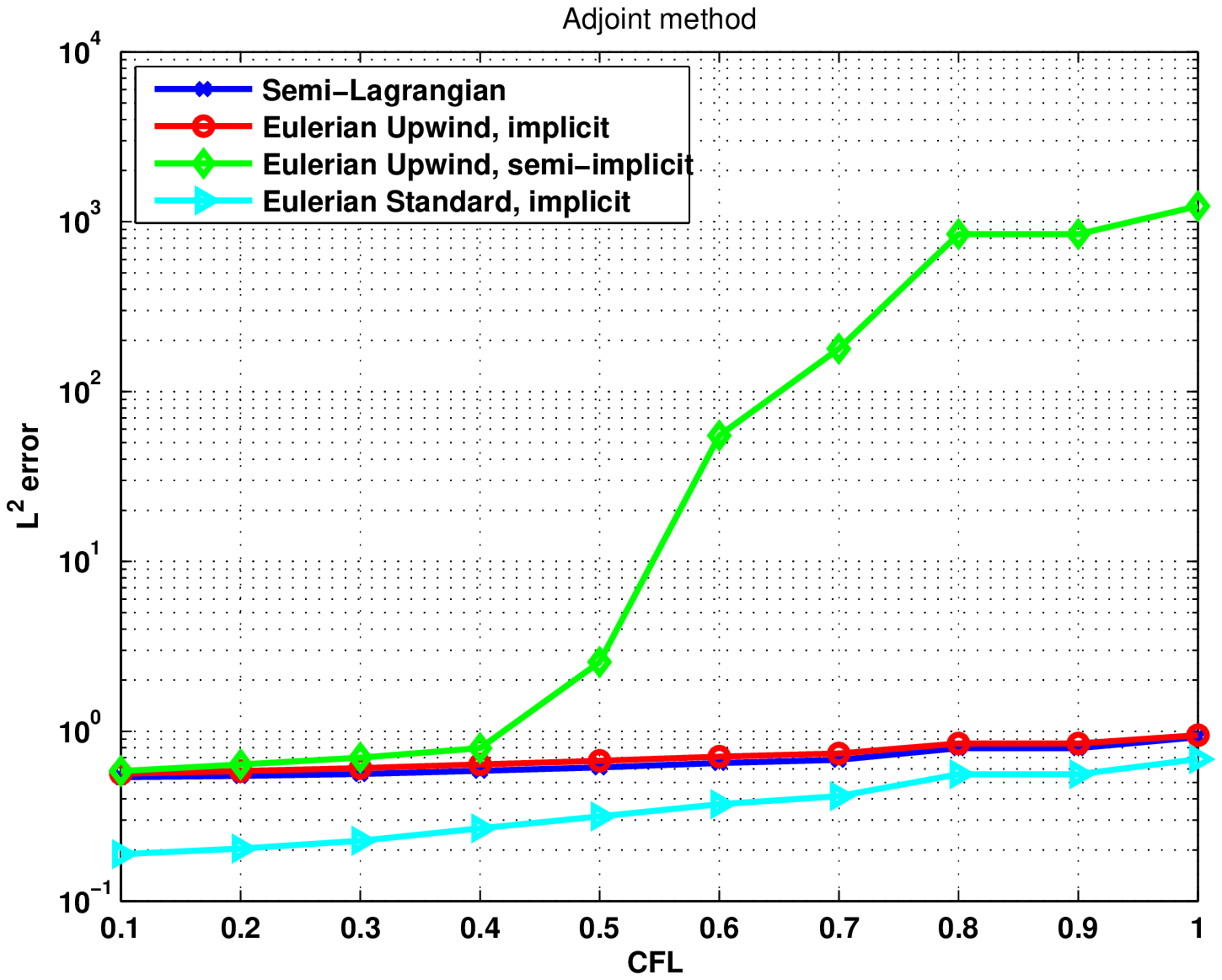}
    \end{center}
     \caption{Experiment III: $L^2$-error for different CFL-numbers at $t=0.5$ for evolution in
     fixed triangulation with mesh size $h_m=0.11$ with
     $\varepsilon=10^{-9}$.}
    \label{fig:stab2}
\end{figure}

\textbf{Experiment III.} Next, we examine the stability properties of the three different schemes for
dominating convection, meaning that we choose  
$\varepsilon$ in problem \eqref{eq:problem1} very small. While for
$\varepsilon=1$ all three schemes are stable for
small and larger CFL-numbers (see fig. (\ref{fig:convplot1}) and
(\ref{fig:convplot2})) we encounter this for small $\varepsilon$ only for the
semi-Lagrangian and the implicit Eulerian scheme (see fig. (\ref{fig:stab1}) and
(\ref{fig:stab2})).

As expected in the convection dominated case, one has to use either the
semi-Lagrangian or fully implicit Eulerian scheme. 


\textbf{Experiment IV.} Another question that we address concerns the stabilization
for the Eulerian schemes. For the implicit timestepping schemes we need to solve in
each time step a stationary convection-diffusion problem.
It is well-known that in the scalar case non-stabilized methods would produce highly
oscillating solutions for convection-diffusion problems with dominating
convection. To study this issue, we consider the stationary 
convection diffusion problem for $1$-forms and $\Omega \subset \mathbb R^2$:
\begin{gather}
    \label{eq:problem2}
    \begin{array}[c]{rcll}
	\ubf + \varepsilon \curl \text{curl} \ubf + \grad({\boldsymbol
	\beta} \cdot \ubf) - \Rbf \boldsymbol \beta \text{curl} \ubf &=& \fbf & \text{in
	}\Omega=[-1,1]^2 \;,\\
	\ubf\cdot\nbf &=& 0 & \text{on }\partial\Omega\;,\\
	\ubf(0) &=& \ubf_{0}\;.
    \end{array}
\end{gather}
We could either use the standard approximation \eqref{eq:limit_standard} or
the upwind interpolated discrete Lie-derivative \eqref{eq:interpol_ufd} in the
discretization of \eqref{eq:problem2}. Again, we choose
\begin{eqnarray*}
  \ubf = 
    \begin{pmatrix} 
	\sin(\pi x)\sin(\pi y) \\  (1-x^2) (1-y^2) 
    \end{pmatrix}.
\end{eqnarray*}
and the divergence free velocity \eqref{eq:velo1}.

\begin{figure}
    \begin{center}
	\includegraphics[width=0.75\linewidth]{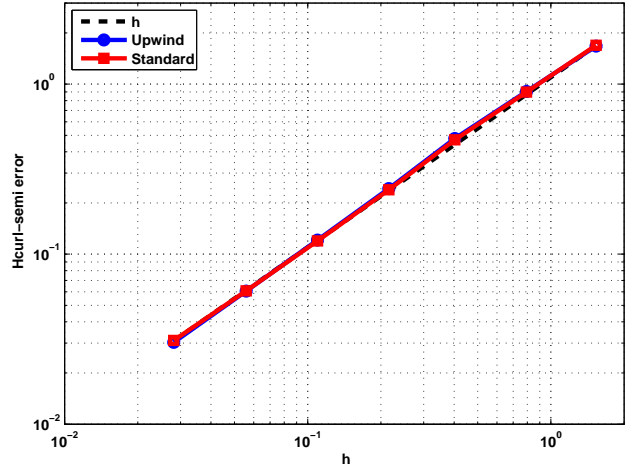}
    \end{center}
    \caption{Experiment IV: $H(\curl)$-error for stationary problem with $\varepsilon=1$.}
    \label{fig:convplot5}
\end{figure}

\begin{figure}
    \begin{center}
	\includegraphics[width=0.75\linewidth]{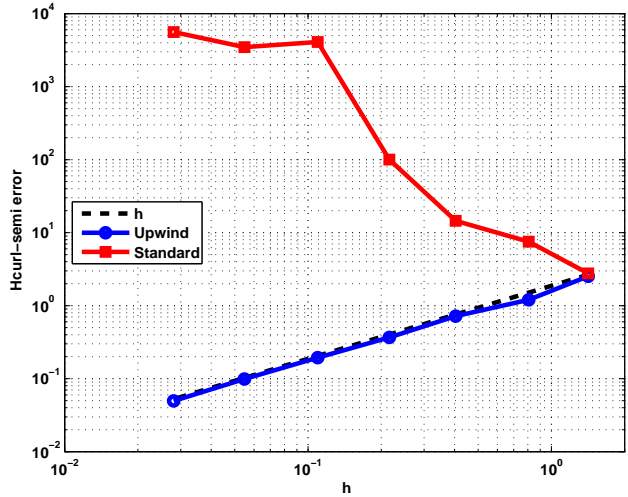}
    \end{center}
    \caption{Experiment IV: $H(\curl)$-error for stationary problem with
   $\varepsilon=10^{-5}$.}
    \label{fig:convplot6}
\end{figure}

Fig.~\ref{fig:convplot5} shows that in the balanced case ($\varepsilon=1$)
both the standard and the stabilizing upwind scheme yield meaningful
solutions.
For small $\varepsilon$ the standard method does not seem to converge while
the upwind scheme converges (see Fig.~\ref{fig:convplot6}).


\end{document}